\documentclass[reqno,11pt]{amsart}
\usepackage{amssymb}
\usepackage{amsmath}
\usepackage{latexsym}
\usepackage{amsfonts}
\newtheorem{Pa}{Paper}[section]
\newtheorem{theorem}[Pa]{{\bf Theorem}}
\newtheorem{algorithm}[Pa]{{\bf Algorithm}}
\newtheorem{lemma}[Pa]{{\bf Lemma}}
\newtheorem{definition}[Pa]{{\bf Definition}}
\newtheorem{corollary}[Pa]{{\bf Corollary}}
\newtheorem{Remark}[Pa]{{\bf Remark}}
\newtheorem{proposition}[Pa]{{\bf Proposition}}
\newtheorem{problem}[Pa]{{\bf Problem}}
\newtheorem{Example}[Pa]{{\bf Example}}

\numberwithin{equation}{section}

\newcommand{\bp}{\boldsymbol{\rho}}
\newcommand{\bgam}{\boldsymbol{\gamma}}
\newcommand{\balpha}{\boldsymbol{\alpha}}
\newcommand{\bbeta}{\boldsymbol{\beta}}

\newcommand{\cZ}{{\mathcal {Z}}}
\newcommand{\cX}{{\mathcal {X}}}
\newcommand{\bH}{{\mathbb H}}
\newcommand{\C}{{\mathbb C}}

\newcommand{\R}{{\mathbb R}}

\newcommand{\bbl}{\boldsymbol{\ell}}
\newcommand{\bbr}{\boldsymbol{r}}   
\newcommand{\bl}{\boldsymbol{e_\ell}}
\newcommand{\br}{\boldsymbol{e_r}}

\begin{document}

\title[]{Zeros and factorizations of quaternion polynomials: the algorithmic 
approach}
\author{Vladimir Bolotnikov}
\address{Department of Mathematics,
College of William and Mary, 
Williamsburg, VA 23187-8795, USA} 
\email{vladi@math.wm.edu}
\maketitle
\begin{abstract}
It is known that polynomials over quaternions may have spherical zeros and isolated left and right 
zeros. These zeros along with appropriately defined multiplicities form the zero structure of a 
polynomial. In this paper, we equivalently describe the zero structure of a polynomial in 
terms of its left and right spherical divisors as well as in terms of left and right indecomposable 
divisors. Several algorithms are proposed to find left/right zeros and 
left/right spherical divisors of a quaternion polynomial, to construct a polynomial with prescribed 
zero structure and more generally, to construct the least left/right common multiple of given 
polynomials. Similar questions are briefly discussed in the setting of quaternion formal power 
series. 
\end{abstract}

\section{Introduction}
\setcounter{equation}{0}
Any monic complex polynomial $p$ of positive degree can be represented as
\begin{equation}
p=\prod_{k=1}^m\bp_{\alpha_k}^{n_k}={\bf lcm}(\bp_{\alpha_1}^{n_1},\ldots,\bp_{\alpha_m}^{n_m}),
\quad \bp_{\alpha_k}:=z-\alpha_k; \; \;  \alpha_k\in\C, \; n_k\in\mathbb N,
\label{1.1}   
\end{equation}
where $\alpha_1,\ldots,\alpha_m$ are {\em distinct zeros} of $p$ of respective multiplicities
$n_1,\ldots,n_m$ and where {\bf lcm} means the least common multiple. 
The set $\{(\alpha_k, n_k)\}_{k=1}^m$ may be referred to as 
to the {\em zero structure} of $p$ and it is clear that a monic polynomial is uniquely 
recovered from its zero structure by formulas \eqref{1.1}.  Since the second representation in 
\eqref{1.1} is equivalent to the primary ideal decomposition 
$\langle p\rangle=\bigcap_{k=1}^m\langle\bp_{\alpha_k}^{n_k}\rangle$
for the ideal $\langle p\rangle$ generated by $p$ and since this decomposition is unique, the zero 
structure of $p$ can be alternatively defined as the collection of generators of irreducible 
components  in the primary ideal decomposition of $\langle p\rangle$. 

\smallskip

The case of polynomials over a division ring (rather than a field) is not that simple 
\cite{wed}, \cite{ore}, \cite{gm}. Here we focus on polynomials over 
the skew  field $\bH$ of quaternions  
\begin{equation}
\alpha=x_0+{\bf i}x_1+{\bf j}x_2+{\bf k}x_3 \qquad (x_0,x_1,x_2,x_3\in\mathbb R),
\label{1.3}   
\end{equation}
where the imaginary units ${\bf i}, {\bf j}, {\bf k}$ commute with $\R$ and satisfy
${\bf i}^2={\bf j}^2={\bf k}^2={\bf ijk}=-1$. 
For $\alpha\in\bH$ of the form \eqref{1.3},
its real and imaginary parts, the quaternion conjugate and the absolute value
are defined as ${\rm Re}(\alpha)=x_0$, ${\rm Im}(\alpha)={\bf i}x_1+{\bf j}x_2+{\bf k}x_3$,
$\overline \alpha={\rm Re}(\alpha)-{\rm Im}(\alpha)$ and
$|\alpha|^2=\alpha\overline{\alpha}=|{\rm Re}(\alpha)|^2+|{\rm Im}(\alpha)|^2$,
respectively. 

\smallskip
\noindent
Two quaternions $\alpha$ and $\beta$ are called {\em equivalent} (conjugate 
to each other) if $\alpha=h^{-1}\beta h$ for some nonzero $h\in\mathbb H$; in notation,
$\alpha\sim\beta$. It turns out (see e.g., \cite{brenner}) that
\begin{equation}
\alpha\sim\beta\quad\mbox{if and only if}\quad {\rm Re}(\alpha) ={\rm Re}(\beta) \;
\mbox{and} \; |\alpha|=|\beta|.
\label{1.4}
\end{equation}
Hence, the {\em conjugacy class} of a given $\alpha\in\mathbb H$ form a $2$-sphere (of radius
$|{\rm Im}(\alpha)|$ around ${\rm Re}(\alpha)$) which will be denoted  by $[\alpha]$.
A finite ordered collection ${\balpha}=(\alpha_1,\ldots,\alpha_k)$ will
be called a {\em spherical
chain} (of the length $k$) if
\begin{equation}
\alpha_1\sim \alpha_2\sim\ldots\sim\alpha_k\quad\mbox{and}\quad \alpha_{j+1}\neq
\overline{\alpha}_j\quad\mbox{for}\quad j=1,\ldots,k-1.
\label{1.8u}
\end{equation}
We let $\bH[z]$ be the ring of polynomials in one formal variable $z$ which commutes with quaternionic 
coefficients. Multiplication in $\bH[z]$ (as well as in $\bH$) is not commutative; however,  
any (left or right) ideal in $\bH[z]$ is principal. We will use notation
$$
\langle p\rangle_{\bf r}:=\left\{pq: \; q\in\bH[z]\right\}\quad\mbox{and}\quad
\langle p\rangle_{\boldsymbol\ell}:=\left\{q p: \; q\in\bH[z]\right\}
$$
for respectively the right and the left ideal generated by $p$. The subscript will be
dropped if the ideal is two-sided; it is not hard to show that any
two-sided ideal in $\bH[z]$ is generated by a polynomial with real coefficients.

\smallskip

For any non-constant polynomial $f\in\bH[z]$, the
sets $\cZ_{\boldsymbol\ell}(f)$ and $\cZ_{\bf r}(f)$ of its left and right zeros
(see Section 2.1 for definitions) are non-empty and are are contained in a finite union of
distinct conjugacy classes:
$$
\cZ(f):=\cZ_{\boldsymbol\ell}(f)\bigcup \cZ_{\bf r}(f)\subset\bigcup V_i.
$$
Moreover, each conjugacy class $V_i$  either contains exactly one left and one right zero
of $f$ or $V_i\subset \cZ_{\boldsymbol\ell}(f)\bigcap \cZ_{\bf r}(f)$; in the latter case,
we say that $V_i$ is the {\em spherical zero} of $f$. Observe that the set of all polynomials 
$f\in\bH[z]$ having the spherical zero $V$ is the two-sided ideal generated by the real polynomial
\begin{equation}
\cX_{V}(z)=(z-\alpha)(z-\overline{\alpha})=z^2-2z\cdot {\rm Re}(\alpha)+|\alpha|^2,\quad \alpha\in V
\label{2.5}
\end{equation}
(the {\em characteristic polynomial} of $V$); it follows from characterization \eqref{1.4}
that $\alpha$ in \eqref{2.5} can be replaced by any other element in $V$.

\smallskip

Several algorithms for finding left/right roots (by means of complex root-finding for a polynomial 
with real coefficients
are known \cite{niven}, 
\cite{ss}, \cite{opj}, \cite{kal}. Yet another algorithm of this type (Algorithm \ref{A:3.4} below) 
is our contribution to the topic. 

\smallskip

By the division algorithm, any non-constant monic $f\in\bH[z]$ can be factored as
\begin{equation}
f(z)=(z-\gamma_1)(z-\gamma_2)\cdots(z-\gamma_N),\quad \gamma_1,\ldots,\gamma_N\in\bH,
\label{1.10}
\end{equation}
where $\gamma_1\in \cZ_{\boldsymbol\ell}(f)$ and $\gamma_N\in\cZ_{\bf r}(f)$. Although
it may happen that none of $\gamma_2,\ldots,\gamma_{N-1}$ belongs to $\cZ(f)$, it is still true that
$\cZ(f)\subset \cup_{k=1}^N[\gamma_k]$.

\smallskip

The {\em zero structure} of quaternion polynomials is a delicate issue (for related results
based on various notions of zero multiplicities, we refer to \cite{ps}, \cite{gs}, \cite{gs1}). Here 
we propose to characterize the left (right) zero structure of a polynomial
in terms of its left (right) {\em spherical divisors}. To introduce these divisors,
we first recall several definitions. 

\smallskip

The {\em least right common multiple} $h={\bf lrcm}(f,g)$ of $f,g\in\bH[z]$ is defined as a (unique)
monic polynomial such that $\langle h\rangle_{\bf r}=\langle f\rangle_{\bf r}\cap \langle g\rangle_{\bf
r}$.  The least left common multiple ${\bf llcm}(f,g)$ is defined as a (unique)
monic polynomial generating the left ideal $\langle f\rangle_{\boldsymbol\ell}\cap  \langle g\rangle_{\boldsymbol\ell}$. 

\smallskip

We will say that a (monic) polynomial $f$ is {\em indecomposable} if it cannot be represented as the 
least right (left) common multiple of its proper left (right) divisors. As we will see in Section 5,
a polynomial $f$ is indecomposable if and only if it admits a {\em unique} factorization \eqref{1.10}
into the product of linear factors and the elements $(\gamma_1,\ldots,\gamma_N)$ from this 
factorization form a spherical chain. In other words, there is a one-to-one correspondence
\begin{equation}
\balpha=(\alpha_{1},\ldots,\alpha_{k}) \mapsto P_{\balpha}:=\bp_{\alpha_{1}}\ldots
\bp_{\alpha_{k}}
\label{3.4u}
\end{equation}
between spherical chains and  monic indecomposable polynomials. Here and in what follows 
we use notation $\bp_\alpha(z):=z-\alpha$ for a fixed $\alpha\in\bH$.
\begin{theorem}
Given $f\in\bH[z]$, let $V_1,\ldots,V_m$ be distinct conjugacy classes containing zeros of $f$.
Then there exist unique (monic) polynomials $D^f_{{\boldsymbol \ell},
V_i}$ and $D^f_{{\boldsymbol r}, V_i}$ ($i=1,\ldots,m$) such that
\begin{align}
f=D^f_{{\boldsymbol \ell}, V_i}\cdot h_i\quad\mbox{so that}\quad
\cZ(D^f_{{\boldsymbol \ell}, V_i})\subset V_i,
\quad \cZ(h_i)\bigcap V_i=\emptyset,\label{1.6}\\
f=g_i\cdot D^f_{{\boldsymbol r}, V_i}\quad\mbox{so that}\quad \cZ(D^f_{{\boldsymbol r},
V_i})\subset V_i,\quad \cZ(g_j)\bigcap V_i=\emptyset.\label{1.7}
\end{align}
Furthermore, $f={\bf lrcm}(D^f_{{\boldsymbol \ell}, V_1},\ldots,D^f_{{\boldsymbol \ell}, V_m})={\bf
llcm}(D^f_{{\boldsymbol r}, V_j},\ldots,D^f_{{\boldsymbol r}, V_m})$ or equivalently,
\begin{equation}
\langle f\rangle_{\bf r} =\bigcap_{i=1}^m\langle D^f_{\bbl,V_i}\rangle_{\bf r},
\quad \langle f\rangle_{\boldsymbol\ell}=\bigcap_{i=1}^m\langle
D^f_{\bbr,V_i}\rangle_{\boldsymbol\ell}.
\label{1.8}
\end{equation}
Finally, there exist (unique) integers $\kappa_i\ge 0$ and spherical chains 
$\balpha_i=(\alpha_{i,1},\ldots,\alpha_{i,k_i})$ and 
$\widetilde\balpha_i=(\widetilde\alpha_{i,1},\ldots,\widetilde\alpha_{i,k_i})$ in $V_i$ such that 
\begin{equation}
D^f_{{\boldsymbol \ell}, V_i}=\cX_{V_i}^{\kappa_i}\cdot P_{\balpha_i}\quad\mbox{and}\quad
D^f_{{\boldsymbol r}, V_j}=P_{\widetilde\balpha_i}\cdot \cX_{V_i}^{\kappa_i} \quad\mbox{for}\quad
i=1,\ldots,m,
\label{1.8a}
\end{equation}
where $P_{\balpha_i}$ and $P_{\widetilde\balpha_i}$ are indecomposable polynomials defined as in \eqref{3.4u}.
Consequently, 
\begin{equation}
f=\cX_{V_1}^{\kappa_1}\cdots \cX_{V_m}^{\kappa_m}\cdot {\bf lrcm}(P_{\balpha_1},\ldots,P_{\balpha_m})=
\cX_{V_1}^{\kappa_1}\cdots \cX_{V_m}^{\kappa_m}\cdot {\bf llcm}(P_{\widetilde\balpha_1},\ldots,P_{\widetilde\balpha_m}).
\label{1.8b}  
\end{equation}
\label{T:1.1} 
\end{theorem} 
We will refer to $D^f_{\bbl,V_j}$ and $D^f_{\bbr,V_j}$ as to {\em left and right spherical divisors}
of $f$. Each left (right) spherical divisor contains all information about zeros of $f$ within the  
corresponding conjugacy  class and thus, we may define the left (right) zero structure of $f$
as the collection of its left (right) spherical divisors. The integer $\kappa_i$ is called  
the {\em spherical zero multiplicity} of the conjugacy class $V_i$. The collections 
$\{(\kappa_i, \balpha_i)\}_{i=1}^m$ and $\{(\kappa_i, \widetilde\balpha_i)\}_{i=1}^m$ 
of spherical multiplicities and spherical chains from representations \eqref{1.8b} also can 
serve as definitions of left and right zero structures of a given polynomial $f$.
A result related to Theorem \ref{T:1.1} is the following version of the Primary Ideal Decomposition Theorem.
\begin{theorem}
For any non-constant polynomial $f\in\bH[z]$, there exist left (right) relatively prime
indecomposable polynomials $p_1,\ldots,p_M$ (resp., $\widetilde{p}_1,\ldots,\widetilde{p}_M$)
such that $f={\bf lcrm}(p_1,\ldots,p_M)={\bf llcm}(\widetilde{p}_1,\ldots,\widetilde{p}_M)$.
Equivalently,   
\begin{equation}
\langle f\rangle_{\bf r} =\bigcap_{i=1}^M\langle p_i\rangle_{\bf r}\quad\mbox{and}\quad
\langle f\rangle_{\boldsymbol\ell}=\bigcap_{i=1}^M\langle
\widetilde{p}_i\rangle_{\boldsymbol\ell}.
\label{1.9}   
\end{equation}
\label{T:1.2} 
\end{theorem}
In contrast to \eqref{1.8}, all ideals in \eqref{1.9} are irreducible, although the polynomials $p_i$ and   
$\widetilde{p}_i$ in \eqref{1.9} are not determined from $f$ uniquely. However, the 
nonuniqueness can be described explicitly (see Theorem \ref{T:5.6} below).

\smallskip

The objective of this paper is to establish explicit connections between representations
\eqref{1.10}, \eqref{1.8}, \eqref{1.8b}, \eqref{1.9} and in particular, to construct explicitly          
a polynomial with the prescribed spherical divisors. One can see from \eqref{1.8b}, that the non-trivial part here 
is to construct the {\bf lrcm} and {\bf llcm} of $m$ indecomposable polynomials with zeros in $m$ distinct
non-real conjugacy classes. The construction has been known for the case where all (but at most one) 
these polynomials are linear (see e.g., \cite{cm}). The algorithm settling the general case is our 
next contribution to the topic. 

\smallskip

The outline of the paper is as follows. After recalling some known results 
on quaternion polynomials in Section 2, we present the root-finding Algorithm \ref{A:3.4} (Section 3)
and Algorithm \ref{A:4.3} which produces the spherical divisor of a given polynomial associated 
with a given conjugacy class.
The proof of Theorem \ref{T:1.1} is given in Section 4.1, and two algorithms relating left
and right spherical divisors corresponding to the same conjugacy class are given in Section 4.2.
In Section 5, we give several characterizations of indecomposable polynomials and prove 
Theorem \ref{T:1.2}. The construction of the 
least common multiple of several indecomposable polynomials is given in 
Lemma \ref{L:5.2} (the case where all polynomials have zeros in the same conjugacy class),
Algorithm \ref{A:6.6} (where the polynomials  have zeros in distinct conjugacy classes).
Finally, Algorithm \ref{A:8.1} produces a polynomial with prescribed spherical divisors 
while the least common multiple of given polynomials is the outcome of Algorithm \ref{A:8.3}.
In the concluding Section 7, we discuss similar questions in the framework of formal power series
over $\bH$.

\section{Background: quaternion polynomials and their zeros}
\setcounter{equation}{0}
A straightforward computation verifies that for any $\alpha\in\bH$ and
$f\in\bH[z]$,
\begin{equation}
f(z)=f^{\bl}(\alpha)+(z-\alpha)\cdot(L_\alpha f)(z)=f^{\br}(\alpha)+(R_\alpha 
f)(z)\cdot(z-\alpha),
\label{2.1}
\end{equation}
where $f^{\bl}(\alpha)$ and $f^{\br}(\alpha)$ are respectively, left and right evaluation of
$f$ at $\alpha$ given by
\begin{equation}
f^{\bl}(\alpha)=\sum_{k=0}^m\alpha^k f_k\quad\mbox{and}\quad
f^{\br}(\alpha)=\sum_{k=0}^m f_k\alpha^k\quad\mbox{if}\quad f(z)=\sum_{j=0}^m z^j f_j,
\label{2.2}
\end{equation}
and where $L_\alpha f$ and $R_\alpha f$ are polynomials of degree $m-1$ given by
\begin{equation}
(L_\alpha f)(z)=\sum_{k=0}^{m-1}\bigg(\sum_{j=0}^{m-k-1}
\alpha^jf_{k+j+1}\bigg)z^k,\quad
(R_\alpha f)(z)=\sum_{k=0}^{m-1}\bigg(\sum_{j=0}^{m-k-1}
f_{k+j+1}\alpha^j\bigg)z^k.
\label{2.3}
\end{equation}
Interpreting $\bH[z]$ as a vector space over $\bH$ we observe that
the mappings $f\mapsto L_\alpha f$ and $f\mapsto R_\alpha f$
define respectively the right linear operator $L_\alpha$ and the left
linear operator $R_\alpha$ (called in analogy to the complex case, the left and the
right backward shift, respectively) acting on $\bH[z]$.

\medskip
\noindent
{\bf 2.1. Left and right zeros.} An element $\alpha\in\bH$ is called a {\em left (right) zero} 
of $f\in\bH[z]$ if $f^{\bl}(\alpha)=0$ (respectively, $f^{\br}(\alpha)=0$).
We denote by $\cZ_{\boldsymbol\ell}(f)$ and $\cZ_{{\boldsymbol r}}(f)$ the sets of  all left and all right zeros of
$f$, respectively, and we let $\cZ(f):=\cZ_{\boldsymbol\ell}(f)\cup \cZ_{\bf r}(f)$. The next two equivalences follow 
from representations \eqref{2.1}:
\begin{equation}
f^{\bl}(\alpha)=0 \; \Leftrightarrow \; f\in \langle \bp_\alpha\rangle_{\bf r},\qquad
f^{\br}(\alpha)=0 \; \Leftrightarrow \; f\in \langle \bp_\alpha\rangle_{\boldsymbol\ell}.\label{2.6}
\end{equation}
Combining equivalences \eqref{2.6} gives the following result (see \cite[Lemma 2.6]{bol1} for 
details).
\begin{lemma}
Let $\alpha,\beta\in\bH$ be two distinct conjugates:
$\beta\in[\alpha]\backslash\{\alpha\}$. Then 
\begin{equation}
\langle \bp_\alpha\rangle_{\bf r}\cap \langle \bp_\beta\rangle_{\bf r}=
\langle \bp_\alpha\rangle_{\boldsymbol\ell}\cap \langle \bp_\beta\rangle_{\boldsymbol\ell}=
\langle \cX_{[\alpha]}\rangle.
\label{2.8}
\end{equation}
{\rm In case the non-real conjugacy class (the $2$-sphere) $V$ is a subset of 
$\cZ_{\boldsymbol\ell}(f)\cap \cZ_{\boldsymbol r}(f)$, we will say that $V$ is a 
{\em spherical zero} of $f$.}
\label{L:2.6}
\end{lemma}

\smallskip
\noindent
{\bf 2.2. Polynomial conjugation.} The quaternionic conjugation $\alpha\mapsto 
\overline{\alpha}$ on $\bH$ can be extended to the anti-linear involution $f\mapsto f^\sharp$ on $\bH[z]$ by letting
\begin{equation}
f^\sharp(z)=\sum_{j=0}^m z^j, \overline{f}_j\quad\mbox{if}\quad f(z)=\sum_{j=0}^m z^j f_j.
\label{2.8a}
\end{equation}
It is not hard to verify that
\begin{equation}
ff^\sharp=f^\sharp f,\quad (fg)^\sharp=g^\sharp f^\sharp,\quad (fg)(fg)^\sharp=(ff^\sharp)(gg^\sharp).
\label{2.10}
\end{equation}
One can see from \eqref{2.2} that if $f\in\R[z]$, then $f^{\bl}(\alpha)=f^{\br}(\alpha)$ for every $\alpha\in\bH$;
in particular, if $f\in\R[z]$, then $\cZ(f)=\cZ_{\boldsymbol\ell}(f)=\cZ_{\bf r}(f)$. Moreover, if
$f\in\R[z]$, then for each $\alpha\in\bH$ and $h\neq 0$,
we have $f(h^{-1}\alpha h)=h^{-1}f(\alpha)h$ so that
$\cZ(f)$ contains, along with each $\alpha$, the whole conjugacy class $[\alpha]$.
Since for any $f\in\bH[z]$, the polynomial $ff^\sharp$ is real, the set $\cZ(ff^\sharp)$ is the
union of finitely many conjugacy classes.  The following result is essentially due to I. Niven 
\cite{niven}.
\begin{theorem}
Let $\deg (f)\ge 1$ and let $\cZ(ff^\sharp)=\bigcup V_i$ be
the union of distinct conjugacy classes.
Then $\cZ_{\boldsymbol\ell}(f)\bigcup \cZ_{\bf r}(f)\subset\cZ(ff^\sharp)$ and each
conjugacy class $V_i$ either contains exactly one left and one right zero of $f$ or
$V_i\in \cZ_{\boldsymbol\ell}(f)\bigcap \cZ_{\bf r}(f)$.
\label{T:2.4} 
\end{theorem}
\noindent
{\bf 2.3. Zero multiplicities:} If we denote by $f^{(k)}$ the $k$-th formal derivative of
$f\in\bH[z]$, then a straightforward verification shows that for any fixed $\alpha\in\bH$,
\begin{equation}
f=\sum_{k=0}^{\deg (f)} \bp_\alpha^k \frac{(f^{(k)})^{\bl}(\alpha)}{k!}=\sum_{k=0}^{\deg
(f)}
\frac{(f^{(k)})^{\br}(\alpha)}{k!}\bp_\alpha^k,\quad \bp_\alpha(z):=z-\alpha.
\label{2.4}
\end{equation}
\begin{definition}
{\rm Let us say that $\alpha\in\bH$ is 
zero of $f\in\bH[z]$ of {\em left zero multiplicity} $m_{\boldsymbol\ell}(\alpha;f)=k$ if 
$f=\bp_\alpha^k h$ for some $h\in\bH[z]$ with $h^{\bl}(\alpha)\neq 0$.}  
\label{D:2.2}
\end{definition}
It follows from \eqref{2.4} that
$f\in \langle \bp_\alpha^k\rangle_{\bf r}\; $ if and only if  $(f^{(j)})^{\bl}(\alpha)=0$ 
for $j=0,\ldots,k-1$; therefore, $k=m_{\boldsymbol\ell}(\alpha;f)$ can be alternatively defined
as the least nonnegative integer such that $(f^{(k)})^{\bl}(\alpha)\neq 0$.  

\smallskip
\noindent
The {\em right zero multiplicity} $m_{\bf r}(\alpha;f)$ is defined as the integer $k$ such that
$f=h\bp_\alpha^k$ for some $h\in\bH[z]$ with $h^{\br}(\alpha)\neq 0$, or equivalently, as the least 
nonnegative integer such that $(f^{(k)})^{\br}(\alpha)\neq 0$.
\begin{definition}
{\rm Let us say that a conjugacy class $V\subset\bH$ is the {\em spherical zero} of $f\in\bH[z]$ of 
the {\em spherical multiplicity} $\kappa=m_s(V;f)$ if 
\begin{equation}
f(z)=\cX_V^\kappa(z) g(z)=g(z)\cX_V^\kappa(z)   
\label{2.16}
\end{equation}
for some $g\in\bH[z]$ vanishing at at most one point in $V$. Equivalently, $\kappa=m_s(V;f)$ 
is the integer such that $f\in\langle\cX_V^\kappa\rangle \backslash \langle\cX_V^{\kappa+1}\rangle$.}
\label{D:2.3}
\end{definition}
The local and spherical zero multiplicities are related as follows: for $f\in\bH[z]$ and a conjugacy class 
$V\subset\bH$,
\begin{equation}
m_s(V;f)=\min_{\gamma\in V} \, \{m_{\boldsymbol\ell}(\gamma;f)\}=\min_{\gamma\in V} \, \{m_{\bf r}(\gamma;
f)\}.
\label{2.15}
\end{equation}
In fact, as it follows by successive applying  Lemma \ref{L:2.6} to $f$ and its formal 
derivatives, both minimums in \eqref{2.15} are attained at all (but at most one) 
elements in $V$. 
\begin{Remark}
Combining \eqref{2.15} with Definition \eqref{D:2.2} of $m_{\boldsymbol\ell}$ and 
$m_{\bf r}$ we conclude that for any two elements $\alpha\neq \beta$ in the conjugacy class
$V\subset \bH$, the spherical zero multiplicity $m_s(V;f)$ equals the least integer $\kappa\ge 0$
such that at least one of the elements $(f^{(\kappa)})^{\bl}(\alpha)$ and $(f^{(\kappa)})^{\bl}(\beta)$
(equivalently, one of $(f^{(\kappa)})^{\br}(\alpha)$ and $(f^{(\kappa)})^{\br}(\beta)$) is non-zero.
\label{R:2.7}
\end{Remark} 
\begin{Example}
Let $\alpha\sim\beta\sim\gamma$ be in the same conjugacy class $V$ and let 
$$
f(z)=(z-\alpha)^n(z-\beta)(z-\gamma)^k.
$$
\begin{enumerate}
\item If $\; \alpha\neq \overline\beta\neq \gamma, \;$ then
$m_{\boldsymbol\ell}(\alpha;f)=n$, $\; m_{\bf r}(\gamma;f)=k$, $\; m_s(V;f)=0$. 
\vspace{1mm}
\item In particular, if $\; \overline\beta\neq \alpha=\gamma, \;$ then 
$\; m_{\boldsymbol\ell}(\alpha;f)=n, \;$ $m_{\bf r}(\alpha;f)=k$,
$\; m_s(V;f)=0$. 
\vspace{1mm}
\item If $\; \gamma=\overline\beta\neq \alpha\neq \beta, \;$ then 
$\; m_{\boldsymbol\ell}(\alpha;f)=n+1, \;$ $m_{\bf r}(\gamma;f)=k$, $\; m_s(V;f)=1$.\vspace{1mm}
\item If $\; \gamma=\alpha=\overline\beta \;$ and $k\le n+1$, then 
$\; m_{\boldsymbol\ell}(\alpha;f)=m_{\bf r}(\alpha;f)=n+1$, $\; m_s(V;f)=k$. 
\end{enumerate}
\label{E:2.4}
\end{Example}
\noindent
{\bf 2.4. Evaluation of products:} By \eqref{2.6}, if $g^{\bl}(\alpha)=0$, then $(gf)^{\bl}(\alpha)=0$
for any $f\in\bH[z]$. On the other hand, since
$(gf)(z)=\sum_{k=0}^n z^k g(z)f_k$, we also have
$$
(gf)^{\bl}(\alpha)=g^{\bl}(\alpha)\sum_{k=0}^n (g^{\bl}(\alpha)^{-1}\alpha
g^{\bl}(\alpha))^kf_k=g^{\bl}(\alpha)f^{\bl}(g^{\bl}(\alpha)^{-1}\alpha g^{\bl}(\alpha)),
$$
provided  $g^{\bl}(\alpha)\neq 0$. Therefore, the left evaluation of the product
of two polynomials is defined by the formula
\begin{equation}
(gf)^{\bl}(\alpha)=\left\{\begin{array}{ccc}
g^{\bl}(\alpha)\cdot f^{\bl}\left(g^{\bl}(\alpha)^{-1}\alpha
g^{\bl}(\alpha)\right)&\mbox{if} &
g^{\bl}(\alpha)\neq 0, \\
0 & \mbox{if} & g^{\bl}(\alpha)= 0.\end{array}\right.
\label{2.17}
\end{equation}
Similarly, $(gf)^{\br}(\alpha)=g^{\br}\left(f^{\br}(\alpha)\alpha
f^{\br}(\alpha)^{-1}\right)\cdot f^{\br}(\alpha)$ if $f^{\br}(\alpha)\neq 0$ and 
$(gf)^{\br}(\alpha)=0$ if $f^{\br}(\alpha)= 0$.
\begin{Remark}
Let $\bgam=(\gamma_1,\ldots,\gamma_N)\subset [\gamma_1]=V$ be a spherical chain.
Then the polynomial 
$f=\bp_{\gamma_1}\bp_{\gamma_2}\cdots\bp_{\gamma_N}$ has a unique left zero $\gamma_1$
and a unique right zero $\gamma_N$. 

\smallskip

Indeed, since $ff^\sharp=\cX_V^N$, it follows by Theorem \ref{T:2.4} that $f$ has no zeros outside 
$V$.
On the other hand, since $(\alpha-\beta)^{-1}\beta(\alpha-\beta)=\overline{\alpha}$ for any two distinct 
$\alpha,\beta\in V$, we conclude from \eqref{2.17} that for every $\gamma\in V\backslash\{\gamma_1\}$,
$$
f^{\bl}(\gamma)=(\gamma-\gamma_{1})(\overline{\gamma}_1-\gamma_{2})(\overline{\gamma}_2-\gamma_3)\cdots   
(\overline{\gamma}_{N-1}-\gamma_{N})\neq 0,
$$
so that $\cZ_{\boldsymbol\ell}(f)=\{\gamma_1\}$. Similarly, one can show that $\cZ_{\bf r}(f)=\{\gamma_N\}$.
\label{R:2.9}
\end{Remark}

\section{A ``new" root-finding algorithm for quaternion polynomials}
\setcounter{equation}{0}

Several known algorithms for finding roots of a quaternion polynomial by means
of complex root-finding (see e.g., \cite{niven}, \cite{ss}, \cite{opj}, \cite{kal}) 
assume that all complex zeros of the real polynomial $ff^\sharp$ are known.   
If $x$ is a real zero of $ff^\sharp$, then $x\in \cZ_{\boldsymbol\ell}(f)\cap  
\cZ_{\boldsymbol r}(f)$, by Theorem \ref{T:2.4}.
If $\alpha,\overline{\alpha}\in\C$ are complex-conjugate roots of $ff^\sharp$ and if
$f^{\bl}(\alpha)=f^{\bl}(\overline{\alpha})=0$,
then $[\alpha]\subset \cZ_{\boldsymbol\ell}(f)\cap \cZ_{\boldsymbol r}(f)$, again by Theorem
\ref{T:2.4}. The remaining case (finding the only left and the only right zero of $f$
in the conjugacy class $[\alpha]$) is the essence of each individual algorithm. Our contribution
here is the following.
\begin{theorem}
Let $f\in\bH[z]$, let $\alpha$ and $\overline\alpha$ be complex roots 
(or any quaternion-conjugate roots) of the real polynomial $ff^\sharp$, and let us assume that 
$f^{\bl}(\alpha)\neq 0$. Then the only left root $\gamma_\ell$ and the only right root $\gamma_r$ of 
$f$ in the conjugacy class $[\alpha]$ are given by
\begin{equation}
\begin{array}{ll}
\gamma_\ell&=(\overline{\alpha}f^{\bl}(\alpha)+\alpha f^{\bl}(\overline{\alpha}))(f^{\bl}(\alpha)+
f^{\bl}(\overline{\alpha}))^{-1},\\ [2mm]
\gamma_r&=(f^{\bl}(\alpha)-f^{\bl}(\overline{\alpha}))^{-1}(\overline{\alpha}f^{\bl}(\alpha)-\alpha
f^{\bl}(\overline{\alpha})).
\end{array}
\label{3.1}
\end{equation}
\label{T:2.1}
\end{theorem}
{\bf Proof:} If $\alpha,\beta,\gamma$ are any three
distinct equivalent quaternions, then
\begin{align}
f^{\bl}(\gamma)=&(\gamma-\beta)(\alpha-\beta)^{-1}f^{\bl}(\alpha)+(\alpha-\gamma)
(\alpha-\beta)^{-1}f^{\bl}(\beta),\label{3.3}\\
f^{\br}(\gamma)=&(\alpha-\beta)^{-1}f^{\bl}(\alpha)\gamma-
\beta(\alpha-\beta)^{-1}f^{\bl}(\alpha)\notag\\
\quad
&+\alpha(\alpha-\beta)^{-1}f^{\bl}(\beta)-(\alpha-\beta)^{-1}f^{\bl}(\beta)\gamma.\label{3.4}
\end{align}
Formula \eqref{3.3} relating evaluations of the same type was established in
\cite{genstr}; for formula \eqref{3.4}, we refer to \cite[Lemma 3.1]{bol1}. Letting $\beta=\overline{\alpha}$
simplifies the latter equalities to
\begin{align}
f^{\bl}(\gamma)=&(\gamma-\overline{\gamma})^{-1}\left[(\gamma-\overline{\alpha})f^{\bl}(\alpha)+
(\gamma-\alpha)f^{\bl}(\overline{\alpha})\right],\label{3.5}\\
f^{\br}(\gamma)=&(\alpha-\overline{\alpha})^{-1}\left[f^{\bl}(\alpha)\gamma-
\overline{\alpha}f^{\bl}(\alpha)
+\alpha f^{\bl}(\overline{\alpha})-f^{\bl}(\overline{\alpha})\gamma\right].\label{3.6}
\end{align}  
Now we will show that $f^{\bl}(\overline{\alpha})\neq \pm f^{\bl}(\alpha)$, so that formulas 
\eqref{3.1} make sense. Indeed, if $f^{\bl}(\overline{\alpha})=f^{\bl}(\alpha)$, then it follows from 
\eqref{3.5} that
$f^{\bl}(\gamma)=f^{\bl}(\alpha)$ for all $\gamma\in[\alpha]$. On the other hand, if
$f^{\bl}(\overline{\alpha})=-f^{\bl}(\alpha)$, then again by \eqref{3.5},
$$
f^{\bl}(\gamma)=(\gamma-\overline{\gamma})^{-1}(\alpha-\overline{\alpha})f^{\bl}(\alpha). 
$$
Neither of the two latter
cases is possible since $f^{\bl}(\alpha)\neq 0$ by the assumption and since $f$ does have a left root
in $[\alpha]$ by Theorem \ref{T:2.4}. Thus $f^{\bl}(\overline{\alpha})\neq \pm f^{\bl}(\alpha)$.

\smallskip

Let $\gamma_\ell$ be a left root of $f$. Then we conclude from \eqref{3.5} that
$$
(\gamma_\ell-\overline{\alpha})f^{\bl}(\alpha)+(\gamma_\ell-\alpha)f^{\bl}(\overline{\alpha})=0,
$$
and solving the latter equation for $\gamma_\ell$ leads us to the first formula in \eqref{3.1}.
Similarly, for the right root $\gamma_r$ we conclude from \eqref{3.6} that 
$$
f^{\bl}(\alpha)\gamma_r- \overline{\alpha}f^{\bl}(\alpha)
+\alpha f^{\bl}(\overline{\alpha})-f^{\bl}(\overline{\alpha})\gamma_r=0,
$$
which being solved for $\gamma_r$ gives the second formula in \eqref{3.1}. \qed
\begin{Remark}
Formulas \eqref{3.1} for unique left and right roots in a given conjugacy class $V$ 
are based on two {\em left} evaluations of $f$ at two quaternion-conjugate points in $V$.
If $f^{\br}(\alpha)\neq 0$, then the formulas for $\gamma_\ell$ and $\gamma_r$ 
in terms of {\em right} evaluations of $f$ are:
\begin{align*}
\gamma_\ell&=(f^{\br}(\alpha)\overline{\alpha}-f^{\br}(\overline{\alpha})\alpha)
(f^{\br}(\alpha)-f^{\br}(\overline{\alpha})^{-1},\\
\gamma_r&=(f^{\br}(\alpha)+f^{\br}(\overline{\alpha}))^{-1}
(f^{\br}(\alpha)\overline{\alpha}+f^{\br}(\overline{\alpha})\alpha).
\end{align*}
The proof is similar to that of Theorem \ref{T:2.1} and will be omitted.
\label{R:3.2}
\end{Remark}
\begin{Remark}
If $f^{\bl}(\overline{\alpha})=0$, then formulas \eqref{3.1}
take the form $\gamma_\ell=\overline{\alpha}$ (as expected) and 
$\gamma_r=(f^{\bl}({\alpha}))^{-1}\overline{\alpha} f^{\bl}({\alpha})$.
Similarly, if $f^{\br}(\overline{\alpha})=0$, then formulas in Remark \ref{R:3.2}
give $\gamma_r=\overline{\alpha}$ and 
$\gamma_\ell=f^{\br}(\alpha)\overline{\alpha}(f^{\br}(\alpha))^{-1}$.
\label{R:2.10}
\end{Remark}
The root-finding algorithm based on the preceeding analysis follows. 
Our contribution to the algorithm is part (3d).
\begin{algorithm}
Given a polynomial $f\in\bH[z]$ of positive degree, 
\begin{enumerate}
\item Find all complex zeros of the real polynomial $ff^\sharp$.\vspace{1mm}
\item Each real zero of $ff^\sharp$ is a zero of $f$.\vspace{1mm}
\item Evaluate $f^{\bl}(\alpha)$ and $f^{\bl}(\overline{\alpha})$ \vspace{1mm}
for each pair $\{\alpha,\overline{\alpha}\}$ of complex-conjugate roots $ff^\sharp$.
\begin{enumerate}
\item If $f^{\bl}(\alpha)=f^{\bl}(\overline{\alpha})=0$, then $[\alpha]$ is the spherical zero of 
$f$.
\item If $f^{\bl}(\alpha)=0\neq f^{\bl}(\overline{\alpha})$, then $\alpha\in\cZ_{\boldsymbol\ell}(f)$
and $(f^{\bl}(\overline{\alpha}))^{-1}\alpha f^{\bl}(\overline{\alpha})\in\cZ_{\bf r}(f)$.
\item If $f^{\bl}(\overline{\alpha})=0\neq f^{\bl}(\alpha)$, then 
$\overline\alpha\in\cZ_{\boldsymbol\ell}(f)$ and $(f^{\bl}(\alpha))^{-1}\overline{\alpha}
f^{\bl}(\alpha)\in\cZ_{\bf r}(f)$.  
\item Otherwise, use formulas \eqref{3.1} to compute $\gamma_\ell\in [\alpha]\cap 
\cZ_{\boldsymbol\ell}(f)$ and  $\gamma_r\in [\alpha]\cap \cZ_{\bf r}(f)$. 
\end{enumerate}
\end{enumerate}
\label{A:3.4}
\end{algorithm}
\begin{Example} To illustrate Algorithm \ref{A:3.4}, take 
$f(z)=z^2-z({\bf j}+2{\bf k})+2{\bf i}$. Then
$$(ff^\sharp)(z)=z^4+5z^2+4=(z^2+1)(z^2+4)
$$
and therefore, all roots of $f$ are contained in the union of two spheres $V_1=[{\bf i}]$ and
$V_2=[2{\bf i}]$. Since
$f^{\bl}({\bf i})=-1+2{\bf i}+2{\bf j}-{\bf k}$ and $f^{\bl}(-{\bf i})=-1+2{\bf i}-2{\bf j}+{\bf k}$,
the formulas \eqref{3.1} give the pair of roots in $V_1$:
$$
\gamma_\ell=(-2{\bf j}-4{\bf k})(-2+4{\bf i})^{-1}={\bf j},\quad
\gamma_r=(-2{\bf k}+4{\bf j})^{-1}(2{\bf i}+4)=0.8{\bf k}-0.6{\bf j}.
$$
Similarly, evaluating  $f^{\bl}({\bf i})=-4+2{\bf i}+4{\bf j}-2{\bf k}$ and $f^{\bl}(-2{\bf i})=-4+2{\bf 
i}-4{\bf j}+2{\bf k}$, we apply formulas \eqref{3.1} to $\alpha=2{\bf i}$ to get another pair of roots
$$
\gamma_\ell=(-8{\bf j}-16{\bf k})(-8+4{\bf i})^{-1}=1.6{\bf j}+1.2{\bf k},\quad
\gamma_r=(8{\bf j}-4{\bf k})^{-1}(8+16{\bf i})=2{\bf k}
$$
in $V_2$. Note that the obtained roots correspond to two different factorizations of $f$:
$$
f(z)=(z-{\bf j})(z-2{\bf k})=\left(z-(1.6{\bf j}+1.2{\bf k})\right)\left(z-(0.8{\bf k}-0.6{\bf j})\right).
$$
\label{E:2.5}
\end{Example}

\section{Factorizations}
\setcounter{equation}{0}

The standard procedure to factorize a (monic) polynomial $f\in\bH[z]$ of degree $N>0$ into
the product $N$ monic linear factors is the following: starting with $Q_0=f$, construct the sequence 
$\{\gamma_j\}\subset\bH$ and the sequence of monic polynomials $\{Q_j\}$ by 
\begin{equation}
\gamma_{j}\in\cZ_{\boldsymbol \ell}(Q_{j-1}),\quad
Q_{j}=L_{\gamma_{j}}Q_{j-1}\quad\mbox{for}
\quad j=1,\ldots,N,
\label{4.1}
\end{equation}
where $L_\alpha$ is the left backward shift operator \eqref{2.3}.
Since $\deg (Q_j)=N-j$ and since $Q_{j-1}^{\bl}(\gamma_j)=0$,
it follows that $Q_{j-1}=\bp_{\gamma_j}Q_j$ and $Q_N\equiv 1$. We now recursively get 
$$
f=Q_0=\bp_{\gamma_1}Q_1=\bp_{\gamma_1}\bp_{\gamma_2}Q_2=\ldots=\bp_{\gamma_1}\bp_{\gamma_2}\cdots 
\bp_{\gamma_N}.
$$
The latter factorization is largely non-unique due to non-unique choices of elements $\gamma_j$ in 
\eqref{4.1}. If we will be picking up the elements $\gamma_j$ from a fixed conjugacy class $V$
for as long as possible, we will recover the zero structure of $f$ within $V$.
Details are furnished in Algorithm \ref{A:4.1} below. Observe that since $ff^\sharp$ is a real polynomial,
$m_s([\alpha];ff^\sharp)=m_{\boldsymbol\ell}(\alpha;ff^\sharp)$.
\begin{algorithm}
Given $f\in\bH[z]$ and a non-real $\alpha\in\cZ(ff^\sharp)$, let $k=m_s([\alpha];ff^\sharp)$. 
\begin{enumerate}
\item Compute recursively $\beta_1,\ldots,\beta_k\in\bH$ and 
$Q_1,\ldots,Q_k\in\bH[z]$ by letting $Q_0=f$ and 
\begin{equation}
\beta_{j+1}=\left\{\begin{array}{ccc}\alpha, &\mbox{if}\; \;  & Q_j^{\bl}(\alpha)=0,\\
(\overline{\alpha}Q_j^{\bl}(\alpha)+\alpha Q_j^{\bl}(\overline{\alpha}))(Q_j^{\bl}(\alpha)+
Q_j^{\bl}(\overline{\alpha}))^{-1},& \mbox{if} \; \;  & Q_j^{\bl}(\alpha)\neq 0,
\end{array}\right. 
\label{4.3}
\end{equation}
\begin{equation}
Q_{j+1}=L_{\beta_{j}}Q_{j}\quad\mbox{for}\quad j=0,\ldots,k-1,
\label{4.4}
\end{equation} 
where $L_{\beta_{j}}$ is the left backward shift operator defined in \eqref{2.3}.
Then
$$
(1) \; \beta_1,\ldots,\beta_k\in[\alpha]; \quad (2) \; \cZ(Q_k)\cap 
[\alpha]=\emptyset;\quad 
(3) \; f=\bp_{\beta_1}\bp_{\beta_2}\cdots \bp_{\beta_k}Q_k.
$$
\item Compute recursively $\widetilde\beta_1,\ldots,\widetilde\beta_k\in\bH$ and
$\widetilde Q_1,\ldots,\widetilde Q_k\in\bH[z]$ by letting $\widetilde Q_0=f$ and
\begin{equation}
\widetilde\beta_{j+1}=\left\{\begin{array}{ccc}\alpha, &\mbox{if} \; \;  & \widetilde 
Q_j^{\bl}(\alpha)=0,\\
(\widetilde{Q}_j^{\bl}(\alpha)-\widetilde{Q}_j^{\bl}(\overline{\alpha}))^{-1}
(\overline{\alpha}\widetilde{Q}_j^{\bl}(\alpha)-{\alpha}\widetilde{Q}_j^{\bl}(\overline{\alpha}))
& \mbox{if}\; \; &\widetilde Q_j^{\bl}(\alpha)\neq 0,
\end{array}\right.
\label{4.6}
\end{equation}
\begin{equation}   
{Q}_{j+1}=R_{\widetilde\beta_{j}}
\widetilde{Q}_{j}\quad\mbox{for}\quad  j=0,\ldots,k-1,
\label{4.7}   
\end{equation}
where $R_{\widetilde\beta_{j}}$ is the right backward shift operator defined in \eqref{2.3}.
Then
$$
(1) \; \widetilde\beta_1,\ldots,\widetilde\beta_k\in[\alpha];\quad (2) \; 
\cZ(\widetilde{Q}_k)\cap [\alpha]=\emptyset; \quad
(3) \; f=\widetilde{Q}_k \bp_{\widetilde\beta_k}\cdots \bp_{\widetilde\beta_2}\bp_{\widetilde\beta_1}.
$$
\end{enumerate}
\label{A:4.1}
\end{algorithm} 
{\bf Proof:} According to Theorem \ref{T:2.1}, the element $\beta_{j+1}$ defined in \eqref{4.3}
is a left zero of the polynomial $Q_j$ and belongs to $[\alpha]$. Therefore, for $Q_{j+1}$ defined 
as in \eqref{4.4},   we have $Q_j=\bp_{\beta_{j+1}}Q_{j+1}$ from which we recursively recover 
$f=Q_0=\bp_{\beta_1}\bp_{\beta_2}\cdots \bp_{\beta_k}Q_k$. By  properties \eqref{2.10}, we conclude from 
the latter representation that $ff^\sharp=\cX_{[\alpha]}^k 
Q_kQ_k^\sharp$ and since $k=m_s([\alpha];ff^\sharp)$, it follows that $Q_k$ has no zeros in 
$[\alpha]$. This completes the proof of the first part of the algorithm. The second part is justified 
in much the same way.\qed

\medskip

Algorithm \ref{A:4.1} produces factorizations requested in \eqref{1.6},
\eqref{1.7} (with $V_j=[\alpha]$, $D_{\boldsymbol\ell,[\alpha]}=\bp_{\beta_1}\cdots 
\bp_{\beta_k}$ and $D_{{\bf r},[\alpha]}=\bp_{\widetilde\beta_k}\cdots \bp_{\widetilde\beta_1}$). 
As was observed in \cite{gs1}, the polynomial $D_{\boldsymbol\ell,[\alpha]}$ (and similarly,
$D_{{\bf r},[\alpha]}$) can be written in a more structured form if some consecutive points obtained via 
\eqref{4.3} are quaternion-conjugates of each other. If $\beta_{j+1}=\overline{\beta}_j$, then 
$\; \bp_{\beta_j}\bp_{\beta_{j+1}}=\bp_{\beta_j}\bp_{\overline{\beta}_{j}}=\cX_{[\alpha]}, \; $
and the latter real polynomial can be commuted through all the factors in $D_{\boldsymbol\ell,[\alpha]}$ 
to the left. Incorporating this observation and Remark \ref{R:4.2} below, we get a more 
efficient modification of Algorithm \ref{A:4.1}.
\begin{Remark}
For $\alpha\in\bH$, a successive application of formulas \eqref{2.3} gives
\begin{equation}
R_\alpha R_{\overline{\alpha}}f=
L_\alpha L_{\overline{\alpha}}f=\sum_{k=0}^{n-2}z^k\sum_{\i=0}^{n-k-2}\left(
\sum_{j=0}^{i}\alpha^j\overline{\alpha}^{i-j}\right)f_{i+k},\quad\mbox{if}\quad f(z)=\sum_{k=0}^n z^kf_k.
\label{4.9} 
\end{equation}
Furthermore, if we define the two-terms recursion
\begin{equation}
r_0=1, \; \; r_1=2{\rm Re}(\alpha), \; \; r_{j+1}=r_j r_1-r_{j-1}|\alpha|^2\quad\mbox{for}\quad 
j=1,2,\ldots,
\label{4.9a}
\end{equation}
an inductive argument shows that $r_k=\sum_{j=0}^{k}\alpha^j\overline{\alpha}^{k-j}$ for all
$k\ge 0$. Therefore, the formula for $L_\alpha L_{\overline{\alpha}}$ in \eqref{4.9}
depends on ${\rm Re}(\alpha)$ and $|\alpha|$ rather than  $\alpha$ itself and therefore, by 
characterization \eqref{1.4}, $L_\alpha L_{\overline{\alpha}}=L_\beta L_{\overline{\beta}}$,
whenever $\alpha\sim\beta$. It thus makes sense to introduce the {\em spherical backward shift}
operator 
\begin{equation}
S_{[\alpha]}f=\sum_{k=0}^{n-2}z^k\sum_{\i=0}^{n-k-2}r_i f_{i+k}\quad\mbox{if}\quad f(z)=\sum_{k=0}^n z^kf_k,
\label{4.9b}  
\end{equation}
where the real numbers $r_i$ are defined in \eqref{4.9a}.
\label{R:4.2}
\end{Remark}
\begin{algorithm}
Given $f\in\bH[z]$ and $\alpha\in\bH$, let $m_s([\alpha];ff^\sharp)=k$.
\begin{itemize}
\item[(1)] Evaluating $f^{(j)}$ at $\alpha$ and $\overline{\alpha}$, find the least integer $\kappa\ge 0$
such that at least one of the elements $(f^{(\kappa)})^{\bl}(\alpha)$ and $(f^{(\kappa)})^{\bl}(\overline\alpha)$
is non-zero.
\item[(2)] Compute $g=S_{[\alpha]}^\kappa f$ by $\kappa$-times application of formula 
\eqref{4.9b}. If $2\kappa=k$, then $g$ has no zeros in $[\alpha]$ and 
$f=\cX^\kappa_{[\alpha]}g=g\cX^\kappa_{[\alpha]}$. 
Otherwise, proceed to $(3a)$ and $(3b)$.
\item[(3a)] Letting $Q_0=g$, use recursive formulas  \eqref{4.3}, \eqref{4.4} $k-2\kappa$ times 
to construct $\alpha_1,\ldots,\alpha_{k-2\kappa}\in\bH$ and polynomials $Q_1,\ldots,Q_{k-2\kappa}=P$. 
\end{itemize}
Then $(\alpha_1,\ldots,\alpha_{k-2\kappa})\subset[a]$ is a spherical chain, $P$ has no zeros in $[\alpha]$, and
\begin{equation}
f=\cX_{[\alpha]}^\kappa\bp_{\alpha_1}\bp_{\alpha_2}\cdots \bp_{\alpha_{k-2\kappa}}P.
\label{4.10}
\end{equation}
\begin{itemize}
\item[(3b)] Letting $\widetilde{Q}_0=g$, use recursive formulas  \eqref{4.6}, \eqref{4.7} $k-2\kappa$ times
to construct $\widetilde{\alpha}_1,\ldots,\widetilde{\alpha}_{k-2\kappa}$ and polynomials 
$\widetilde{Q}_1,\ldots,\widetilde{Q}_{k-2\kappa}=\widetilde{P}$.
\end{itemize}
Then $(\widetilde{\alpha}_1,\ldots,\widetilde{\alpha}_{k-2\kappa})\subset [\alpha]$ is a spherical chain,  
$\widetilde{P}$ has no zeros in $[\alpha]$, and
\begin{equation}
f=\widetilde{P}\bp_{\widetilde{\alpha}_{k-2\kappa}}\cdots \bp_{\widetilde{\alpha}_{2}}\bp_{\widetilde{\alpha}_{1}}
\cX_{[\alpha]}^\kappa.
\label{4.11}
\end{equation}
\label{A:4.3}
\end{algorithm}
{\bf Proof:} The integer $\kappa$ obtained in Step (1) equals  
$m_s([\alpha];f)$ by Remark \ref{R:2.7}. Hence the polynomial $g=S_{[\alpha]}^\kappa f$
satisfies \eqref{2.16}, and $m_s([\alpha];g)=0$. On the other hand, since $ff^\sharp=\cX_{[\alpha]}^{2\kappa}
gg^\sharp$ (by \eqref{2.10}), it follows that $m_s([\alpha];gg^\sharp)=k-2\kappa$ and we may apply Algorithm
\ref{A:4.1} (part 1) to the polynomial $g$ to get $\alpha_1,\ldots,\alpha_{k-2\kappa}\in[\alpha]$ and 
the polynomial $P_{k-2\kappa}$ with no zeros in $[\alpha]$ such that 
$g=\bp_{\alpha_1}\bp_{\alpha_2}\cdots \bp_{\alpha_{k-2\kappa}}P_{k-2\kappa}$.    
If $\alpha_{j+1}=\overline{\alpha}_j$ for some $j$, then $\bp_{\alpha_j}\bp_{\alpha_{j+1}}=\cX_{[\alpha]}$
so that $m_s([\alpha];g)\ge 1$ which is a contradiction. Therefore $\alpha_{j+1}\neq \overline{\alpha}_j$
for all $j$ as desired. Applying the second part of Algorithm
\ref{A:4.1} to the polynomial $g$, we get representation \eqref{4.11}.\qed
\begin{Remark}
The representation \eqref{4.10} is unique in the following sense: {\em if
\begin{equation} 
f=\cX_{[\alpha]}^{\kappa^\prime}\bp_{\alpha^\prime_1}\bp_{\alpha^\prime_2}\cdots \bp_{\alpha^\prime_{s}}G,
\qquad \cZ_{\boldsymbol\ell}(G)\cap[\alpha]=\emptyset,
\label{4.12}
\end{equation}
is another factorization of $f$ with the spherical chain
$(\alpha^\prime_1,\ldots,\alpha^\prime_{s})\in [\alpha]$, then $\kappa^\prime=\kappa$, $s=k-2\kappa$, 
$\alpha_j^\prime=\alpha_j$ for $j=1,\ldots,s$, and $G=P$}. The representation
\eqref{4.11} is unique in a similar sense.
\label{R:4.4}
\end{Remark}
{\bf Proof:} The integers $\kappa$ and $\kappa^\prime$ are both equal to $m_s([\alpha],f)$ and therefore, 
$\kappa^\prime=\kappa$. Then we have from \eqref{4.10} and \eqref{4.12}, 
\begin{equation}
g=\bp_{\alpha_1}\bp_{\alpha_2}\cdots\bp_{\alpha_{k-2\kappa}}P_{k-2\kappa}=\bp_{\alpha^\prime_1}
\bp_{\alpha^\prime_2}\cdots \bp_{\alpha^\prime_{s}}G.
\label{4.13}
\end{equation}
Therefore $gg^\sharp=\cX_{[\alpha]}^{2k-4\kappa}P_{k-2\kappa}P^\sharp_{k-2\kappa}=\cX_{[\alpha]}^{2s}GG^\sharp$
and since $P_{k-2\kappa}$ and $G$ have no zeros in $[\alpha]$, we conclude that $k-2\kappa=s$.
By Remark \ref{R:2.9}, it follows from factorizations \eqref{4.13} that $\alpha_1$ (and also $\alpha^\prime_1$) 
is a unique left zero of $g$ in $V$. Therefore, $\alpha_1=\alpha_1^\prime$ and applying $L_{\alpha_1}$ to equalities 
\eqref{4.13} gives
$$
L_{\alpha_1}g=\bp_{\alpha_2}\cdots \bp_{\alpha_{s}}P=\bp_{\alpha^\prime_2}\cdots\bp_{\alpha^\prime_{s}}G.
$$
Repeating the above argument we subsequently conclude that $\alpha_j=\alpha_{j}^\prime$ for all 
$j=1,\ldots,s$ and then also $P=G$.\qed
\begin{corollary}
Any monic polynomial $f$ with all zeros contained in the non-real conjugacy class $V$
can be (uniquely) factored either as $f=\cX_{[\alpha]}^\kappa$ or as 
\begin{equation}
f=\cX_{V}^\kappa\bp_{\alpha_1}\cdots\bp_{\alpha_{n}}\qquad
(\alpha_j\in V, \; \alpha_{j+1}\neq \overline\alpha_j).
\label{4.14}
\end{equation}
{\rm The statement follows from representation \eqref{4.10} and Remark \ref{4.4}, since the monic 
polynomial $P$ in \eqref{4.10} does not have roots and therefore $P\equiv 1$.} 
\label{C:4.5}
\end{corollary}
\begin{Example}
To illustrate Algorithm \ref{A:4.3}, let us consider the polynomial
\begin{align*}
f(z)=&z^7-(1+{\bf i}+{\bf j}+{\bf k})z^6+(2-{\bf i}+2{\bf j})z^5-(3+{\bf i}+2{\bf j}+2{\bf
k})z^4\\&+(1-2{\bf i}+4{\bf j})z^3-(3-{\bf i}+{\bf j}+{\bf k})z^2+(2{\bf j}-{\bf i})z+{\bf i}-1.
\end{align*}
A straightforward computation shows that
\begin{align*}
(ff^\sharp)(z)=&z^{14}-2z^{13}+8z^{11}+27z^{10}
-30z^{9}+50z^{8}-40z^{7}+55z^{6}-30z^{5}\\
&+36z^{4}-12z^{3}+13z^{2}- 2z + 2=(z^2+1)^6(z^2-2z+2).
\end{align*}
We see that all zeros of $f$ are contained in the conjugacy classes
$V_1=[{\bf i}]$ and $V_2=[1+{\bf i}]$. $V_2$ contains isolated zeros
of multiplicity one, whereas $V_1$ contains zeros of higher multiplicities.
Applying Algorithm \ref{A:4.3} we first evaluate:
$f^{\bl}({\bf i})=f^{\bl}(-{\bf i})=0$, $\, (f^{\prime})^{\bl}({\bf i})=(f^{\prime})^{\bl}(-{\bf
i})=0$ and
\begin{align*}
(f^{\prime\prime})^{\bl}({\bf i})=&42{\bf i}^5-30{\bf i}^4(1+{\bf i}+{\bf j}+{\bf k})
+20{\bf i}^3(2-{\bf i}+2{\bf j})-12{\bf i}^2(3+{\bf i}+2{\bf j}+2{\bf k})\\
&+6{\bf i}(1-2{\bf i}+4{\bf j})-2(3-{\bf i}+{\bf j}+{\bf k})=
-8-8{\bf i}-8{\bf j}-24{\bf k}\neq 0.
\end{align*}
Thus, $m_s([{\bf i}], f)=2$. The recursion \eqref{4.9a} takes the form
$$
r_0=1,\quad r_i=0,\quad r_{j+1}=-r_{j-1}\quad\mbox{for}\quad j=1,2,\ldots
$$
and we subsequently get
\begin{align*}
S_{[{\bf i}]}f&=z^5-(1+{\bf i}+{\bf j}+{\bf k})z^4+(1-{\bf i}+2{\bf j})z^3
-(2+{\bf j}+{\bf k})z^2+(2{\bf j}-{\bf i})z+{\bf i}-1,\\
g&:=S^2_{[{\bf i}]}f=z^3-(1+{\bf i}+{\bf j}+{\bf k})z^2-({\bf i}-2{\bf j})z+{\bf i}-1.
\end{align*}
Since $m_s([{\bf i}], ff^\sharp)=4$ and  $m_s([{\bf i}], f)=2$, we proceed to Step 3. We have
\begin{equation}
g^{\bl}({\bf i})=1+{\bf i}+{\bf j}+3{\bf k}, \quad g^{\bl}({\bf -i})=-1+3{\bf i}+{\bf j}-{\bf k},
\label{4.17}
\end{equation}
and therefore, by part (3a) of Algorithm \ref{A:4.3},
\begin{align*}
\alpha_1=&(-{\bf i}g^{\bl}({\bf i})+{\bf i}g^{\bl}({\bf -i}))(g^{\bl}({\bf i})+g^{\bl}({\bf -i}))^{-1}={\bf k},\\
Q_1=&L_{\alpha_1}g=L_{\bf k}g=z^2-(1+{\bf i}+{\bf j})z+{\bf j}-{\bf k}.
\end{align*}
We next compute
$Q_1^{\bl}({\bf i})=-{\bf i}+{\bf j}-2{\bf k},\quad Q_1^{\bl}(-{\bf i})=-2+{\bf i}+{\bf j}$ and subsequently get
\begin{align*}
\alpha_2=&(-{\bf i}Q_1^{\bl}({\bf i})+{\bf i}Q_1^{\bl}({\bf -i}))(Q_1^{\bl}({\bf i})+Q_1^{\bl}({\bf -i}))^{-1}={\bf j},\\
Q_2=&L_{\alpha_2}Q_1=L_{\bf j}Q_1=z-1-{\bf i}.
\end{align*}
The representation \eqref{4.10} for $f$ takes the form
\begin{equation}
f(z)=(z^2+1)^2(z-{\bf k})(z-{\bf j})P(z),\quad\mbox{where}\quad P(z)=(z-1-{\bf i}).
\label{4.17a}
\end{equation}
Applying part (3b) of Algorithm \ref{A:4.3} and making use of \eqref{4.17}, we get
\begin{align*}
\widetilde{\alpha}_1=&(g^{\bl}({\bf i})-g^{\bl}({\bf -i}))^{-1}(-{\bf i}g^{\bl}({\bf i})-{\bf
i}g^{\bl}({\bf -i}))=\frac{2{\bf i}+{\bf j}-2{\bf k}}{3},\\
\widetilde{Q}_1=&R_{\widetilde{\alpha}_1}g=
z^2-\frac{3+{\bf i}+2{\bf j}+5{\bf k}}{3}z+\frac{-2-2{\bf i}+{\bf j}+3{\bf k}}{3}.
\end{align*}
We then compute
\begin{align*}
\widetilde{Q}_1^{\bl}({\bf i})=&\frac{-4-5{\bf i}+6{\bf j}+{\bf k}}{3},\qquad
\widetilde{Q}_1^{\bl}(-{\bf i})=\frac{-6+{\bf i}-4{\bf j}+5{\bf k}}{3},\\
\widetilde{\alpha}_2=&(\widetilde{Q}_1^{\bl}({\bf i})-\widetilde{Q}_1^{\bl}({\bf -i}))^{-1}(-{\bf
i}\widetilde{Q}_1^{\bl}({\bf i})-{\bf
i}\widetilde{Q}_1^{\bl}({\bf -i}))=\frac{-2{\bf i}+26{\bf
j}+29{\bf k}}{39},\\
\widetilde{Q_2}=&R_{\widetilde{\alpha}_2}\widetilde{Q}_1=z-1+\frac{5{\bf i}+12{\bf k}}{13},
\end{align*}
and representation \eqref{4.11} for $f$ takes the form
\begin{equation}
f(z)=\widetilde{P}(z)\left(z-\frac{-2{\bf i}+26{\bf
j}+29{\bf k}}{39}\right)\left(z-\frac{2{\bf i}+{\bf j}-2{\bf
k}}{3}\right)(z^2+1)^2,
\label{4.17b}
\end{equation}
where $\widetilde{P}(z)=z-1-\frac{5{\bf i}+12{\bf k}}{13}$.
\label{E:4.3}  
\end{Example}
{\bf 4.1. Proof of Theorem \ref{T:1.1}.} Let $V$ be a conjugacy class containing zeros of a 
given $f\in\bH[z]$.
If $V=\{x\}$ where $x$ is a real root of $f$ of multiplicity $k$, then $f=\bp_x^k h=h\bp_x^k$ are 
the factorizations requested in \eqref{1.6}, \eqref{1.7}. In the non-real case, representations 
\eqref{1.6}, \eqref{1.7} are established by Algorithm \ref{A:4.3} with 
\begin{equation}
D^f_{\boldsymbol \ell,[\alpha]}=\cX_{[\alpha]}^\kappa\bp_{\alpha_1}\cdots 
\bp_{\alpha_{k-2\kappa}}\quad\mbox{and}\quad D^f_{{\bf r},[\alpha]}=\cX_{[\alpha]}^\kappa
\bp_{\widetilde{\alpha}_{k-2\kappa}}\cdots\bp_{\widetilde{\alpha}_{1}},
\label{4.15}
\end{equation}
proving, therefore, also formulas \eqref{1.8a}. The uniqueness of factorizations  \eqref{1.6}, \eqref{1.7} 
and \eqref{1.8a} was shown in Remark \ref{4.4} and Corollary \ref{C:4.5}.
Since the real polynomial $ff^\sharp$ has only spherical zeros or isolated real zeros of even multiplicities,
it can be factored as $ff^\sharp=\prod_{j=1}^m\cX_{V_j}^{k_j},$
where we again let $\cX_{[x]}=\bp_{x}^2$ if $x\in\mathbb R$. It follows from Algorithm \ref{A:4.3}
and formulas \eqref{4.15} that $\deg (D^f_{\boldsymbol\ell,V_j})=k_j$ so that $\deg f= 
\sum_{j=1}^m k_j$.
Observe that for any right common multiple
$F$ of $D^f_{\boldsymbol\ell,V_1},\ldots,D^f_{\boldsymbol\ell,V_m}$, the polynomial
$FF^\sharp$ is a common multiple of relatively prime real polynomials
$D^f_{\boldsymbol\ell,V_j}D_{\boldsymbol\ell,V_j}^\sharp=\cX_{V_j}^{k_j}$
($1\le j\le m$) and therefore, $\deg F\ge \sum_{j=1}^m k_j$. Thus, $f$ is a 
a right common multiple of $D^f_{\boldsymbol\ell,V_1},\ldots,D^f_{\boldsymbol\ell,V_m}$ 
(by \eqref{1.6}) of the minimally possible degree. Therefore, $f={\bf lrcm}(D^f_{{\boldsymbol \ell}, 
V_1},\ldots,D^f_{{\boldsymbol \ell}, V_m})$ which proves the first equality in \eqref{1.8}. The second 
equality follows similarly. Equalities \eqref{1.8b} now follow from \eqref{1.8a}.\qed

\medskip
\noindent
{\bf 4.2. The right zero structure versus the left:} If $f\in\bH[z]$ is completely factored
as in \eqref{1.10}, we can construct its spherical divisors using Algorithm \ref{A:4.3}.
If $f$ is given in the form \eqref{4.10}, then its left zero structure is known only within 
the conjugacy class $[\alpha]$. However, this information is sufficient to recover the right zero 
structure of $f$ within $[\alpha]$. We recall the backward shift operators $L_{\alpha}$ and 
$R_{\beta}$ defined in \eqref{2.3}. 
\begin{lemma}
{\rm (1)} If $F\in\bH[z]$ and $\gamma\in\bH$ are such that
$\cZ(F)\cap[\gamma]=\emptyset$, then
$$
\bp_{\gamma}F=Q\bp_{\beta},\quad\mbox{where}\quad
\beta=(F^{\bl}(\overline{\gamma}))^{-1}\gamma F^{\bl}(\overline{\gamma}),
\quad Q=R_{\beta}(\bp_{\gamma}F).
$$
Moreover, $\cZ(Q)\cap[\gamma]=\emptyset$ and $\gamma=F^{\br}(\beta)\beta (F^{\br}(\beta))^{-1}$.

\medskip

{\rm (2)} If $Q\in\bH[z]$ and $\beta\in\bH$ are such that $\cZ(Q)\cap[\beta]=\emptyset$, then
$$
Q\bp_{\beta}=\bp_{\gamma}F,\quad\mbox{where}\quad
\gamma=Q^{\br}(\overline{\beta})\beta (Q^{\br}(\overline{\beta}))^{-1},
\quad F=L_{\gamma}(Q\bp_{\beta}).
$$
Moreover, $\cZ(F)\cap[\beta]=\emptyset$ and $\beta=(Q^{\bl}(\gamma))^{-1}\gamma Q^{\bl}(\gamma)$.
\label{L:4.7}
\end{lemma} 
{\bf Proof:} The polynomial $g=\bp_{\gamma}F$ has a unique
left zero in $[\gamma]$ (which is  $\gamma$). By Remark \ref{R:2.10}, the unique right zero of $g$ in
$[\gamma]$ is given by
$$
\beta=g^{\bl}(\overline{\gamma})^{-1}\gamma g^{\bl}(\overline{\gamma})=
((\overline{\gamma}-\gamma)F^{\bl}(\overline{\gamma}))^{-1}\gamma
(\overline{\gamma}-\gamma)F^{\bl}(\overline{\gamma})
=(F^{\bl}(\overline{\gamma}))^{-1}\gamma F^{\bl}(\overline{\gamma}).
$$
By \eqref{2.1}, $g$ can be factored as $g=Q\bp_{\beta}$, where
$Q=R_{\beta}g$. Since $\beta$ is the only right zero of $g$ 
in $[\gamma]$, it follows that $\cZ(Q)\cap[\gamma]=\emptyset$. 
Finally, evaluating both parts in $\bp_{\gamma}F=Q\bp_{\beta}$
at $z=\beta$ on the right gives
$$
F^{\br}(\beta)\beta-\gamma F^{\br}(\beta)=0
$$
from which we conclude $\gamma=F^{\br}(\beta)\beta (F^{\br}(\beta))^{-1}$.
This completes the proof of the first statement of the lemma. 
The second statement is verified in much the same way.\qed 

\medskip

The next algorithm recovers the right spherical divisor of a given polynomial from the 
given left spherical divisor associated with the same conjugacy class.
\begin{algorithm}
Given $f\in\bH[z]$ in the form 
\begin{equation}
f=D^f_{\boldsymbol \ell,[\alpha]}P,\quad\mbox{where}\quad D^f_{\boldsymbol
\ell,[\alpha]}=\cX_{[\alpha]}^\kappa\bp_{\alpha_1}\bp_{\alpha_2}\cdots \bp_{\alpha_{n}}\quad
(\alpha_{j+1}\neq \overline{\alpha}_j),
\label{4.20}
\end{equation}
let $P_0=P$ and recursively compute 
\begin{equation}
\widetilde{\alpha}_{j+1}=P^{\bl}_j(\overline{\alpha}_{n-j})^{-1}\alpha_{n-j} 
P^{\bl}_j(\overline{\alpha}_{n-j}),\quad P_{j+1}=R_{\widetilde{\alpha}_{j+1}}(\bp_{\alpha_{n-j}}P_j)
\label{4.21}
\end{equation}
for $j=0,\ldots,n-1$. Then $f$ can be represented as 
\begin{equation}
f=\widetilde{P}D^f_{{\bf r},[\alpha]},\quad\mbox{where}\quad 
D^f_{{\bf r},[\alpha]}=\bp_{\widetilde{\alpha}_{n}}\cdots \bp_{\widetilde{\alpha}_{2}}
\bp_{\widetilde{\alpha}_{1}}\cX_{[\alpha]}^\kappa, \quad \widetilde{P}=P_n.
\label{4.22}
\end{equation}
\label{A:4.8}
\end{algorithm}
{\bf Proof:} Based on the first statement in Lemma \ref{L:4.7} and definitions 
\eqref{4.21}, a simple induction argument shows that
$$
\bp_{\alpha_{n-j}}P_j=P_{j+1}\bp_{\widetilde\alpha_{j}}\quad\mbox{and}\quad \cZ(P_{j+1})\cap[\alpha]=\emptyset
\quad\mbox{for}\quad j=0,\ldots,n-1.
$$ 
Therefore, we get recursively
\begin{align*}
\bp_{\alpha_1}\cdots \bp_{\alpha_{n-1}}\bp_{\alpha_{n}}P&
=\bp_{\alpha_1}\cdots\bp_{\alpha_{n-1}} \bp_{\alpha_{n}}P_0
=\bp_{\alpha_1}\cdots\bp_{\alpha_{n-1}}P_{1}\bp_{\widetilde{\alpha}_1}\\
&=\bp_{\alpha_1}\cdots\bp_{\alpha_{n-2}}P_2\bp_{\widetilde{\alpha}_2}\bp_{\widetilde{\alpha}_1}
=\ldots 
=P_{n}\bp_{\widetilde{\alpha}_n}\bp_{\widetilde{\alpha}_{n-1}}\cdots \bp_{\widetilde{\alpha}_1}.
 \end{align*}
Multiplying the last equality by $\cX_{[\alpha]}^\kappa$ (from either side, since $\cX_{[\alpha]}\in\mathbb 
R[z]$) and taking into account \eqref{4.20} we get \eqref{4.22}. Since $m_s([\alpha;f])=\kappa$ by 
\eqref{4.20}, it also follows that $\widetilde{\alpha}_{j+1}\neq \overline{\widetilde\alpha}_j$ for 
$j=0,\ldots,n-1$.\qed
\begin{Example}
Let us consider the polynomial $f$ from Example \ref{E:4.3} and its factorization 
\eqref{4.17a}. Thus, $\alpha_1={\bf k}$, ${\alpha_2}={\bf j}$ and $P(z)=z-1-{\bf i}$.
We now apply Algorithm \ref{A:4.8}:
\begin{align*}
\widetilde{\alpha}_1&=(-{\bf j}-1-{\bf i})^{-1}{\bf j}(-{\bf j}-1-{\bf i})=\frac{2{\bf i}+{\bf j}-2{\bf k}}{3},\\
\bp_{\alpha_2}P_0&=z^2-(1+{\bf i}+{\bf j})z+{\bf j}-{\bf k},\quad P_1=R_{\widetilde{\alpha}_1}(\bp_{\alpha_2}P_0)=
z-1-\frac{{\bf i}+2{\bf j}+2{\bf k}}{3},\\
\widetilde{\alpha}_2&=\left(-{\bf k}-1-\frac{{\bf i}+2{\bf j}+2{\bf k}}{3}\right)^{-1}{\bf k}\left(-{\bf k}-1-\frac{{\bf 
i}+2{\bf j}+2{\bf k}}{3}\right)=\frac{-2{\bf i}+26{\bf j}+29{\bf k}}{39},\\
\bp_{\alpha_1}P_1&=z^2-\frac{3+{\bf i}+2{\bf j}+5{\bf k}}{3}z+\frac{-2-2{\bf i}-{\bf j}+3{\bf k}}{3},\\
P_2&=R_{\widetilde{\alpha}_2}(\bp_{\alpha_1}P_1)=z-1-\frac{5{\bf i}+12{\bf k}}{13},
\end{align*}
and the factorization $f=P_2 \bp_{\widetilde{\alpha}_2}\bp_{\widetilde{\alpha}_1}\cX_{[{\bf i}]}^2$ coincides with
that in \eqref{4.17b}, as expected.
\label{E:4.10}
\end{Example}
We conclude the section with the algorithm recovering $D^f_{\boldsymbol \ell,[\alpha]}$ from $D^f_{{\bf
r},[\alpha]}$; justification is based on the second statement in Lemma \ref{L:4.7} and will be omitted. 
\begin{algorithm}
Given $f\in\bH[z]$ of the form \eqref{4.22},
let $P_0=P$ and recursively compute
\begin{equation}
\alpha_{j+1}=P^{\br}_j(\overline{\widetilde\alpha}_{n-j})^{-1}\widetilde\alpha_{n-j}
P^{\br}_j(\overline{\widetilde\alpha}_{n-j}),\quad P_{j+1}=
L_{\alpha_{j+1}}(P_j\bp_{\widetilde\alpha_{n-j}})
\label{4.23}
\end{equation}
for $j=0,\ldots,n-1$. Then $f$ can be represented as in \eqref{4.20} with $P=P_n$.
\label{A:4.9}
\end{algorithm}

\section{Indecomposable polynomials and irreducible decompositions}
\setcounter{equation}{0}

Let us say that a left (right) ideal in $\bH[z]$ is irreducible if it is not contained 
properly in two distinct ideals of the same type. The generators of irreducible ideals, 
therefore, are the polynomials that cannot be represented as the least right (left)
common multiple of their proper left (right) divisors. In \cite{ore}, 
such polynomials were called {\em indecomposable}. In the next theorem, we collect
a number of equivalent characterizations of indecomposable  polynomials.
\begin{theorem}
Let $f\in\bH[z]$ be a monic polynomial factored as in \eqref{1.10}:
\begin{equation}
f(z)=(z-\gamma_1)(z-\gamma_2)\cdots(z-\gamma_N),\quad \gamma_1,\ldots,\gamma_N\in\bH.
\label{5.1}
\end{equation}
The following are equivalent:
\begin{enumerate}
\item $\bgam=(\gamma_1,\ldots,\gamma_N)$ is a spherical chain.
\item $\gamma_1$ is the only left zero of $f$.
\item $\gamma_N$ is the only right zero of $f$.
\item \eqref{5.1} is a unique factorization of $f$ into the product of linear factors.
\item The ideal $\langle f\rangle_{\bf r}$ is irreducible.
\item The ideal $\langle f\rangle_{\boldsymbol\ell}$ is irreducible.
\end{enumerate}
\label{T:5.1} 
\end{theorem}
{\bf Proof:} The implication $(1)\Rightarrow(2)$ was verified in Remark \ref{R:2.9}. 
Let us assume that (1) is not in force, i.e., that either $\gamma_1\not\sim\gamma_j$
or $\gamma_{j+1}=\overline{\gamma}_j$ for some $j\in\{1,\ldots,N\}$. In the first case, the
conjugacy class $[\gamma_j]$ contains a left zero $\alpha$ of $f$ different from $\gamma_1$; 
in the second case, $\bp_{\gamma_j}\bp_{\gamma_{j+1}}=\bp_{\gamma_j}\bp_{\overline{\gamma}_{j}}=\cX_V$ so that
$f\in\langle\cX_V\rangle$ and therefore $f$ has infinitely many left zeros. This completes the proof of 
$(1)\Leftrightarrow(2)$.

\smallskip

Let us assume that (1) holds and let 
$f(z)=(z-\gamma_1^\prime)(z-\gamma_2^\prime)\cdots(z-\gamma_N^\prime)$ be another 
factorization of $f$ into the product of linear factors. Since both $\gamma_1$ and $\gamma_1^\prime$ 
are left zeros of $f$ and  since $\gamma_1$ is the only left zero of $f$ by $(1)\Rightarrow(2)$, it 
follows that $\gamma_1=\gamma_1^\prime$. We then consider the equality 
$$
L_{\gamma_1}f=\bp_{\gamma_2}\cdots \bp_{\gamma_N}=\bp_{\gamma_2^\prime}\cdots \bp_{\gamma^\prime_N}
$$
and since still $\gamma_{j+1}\neq\overline{\gamma}_j$ for $j=2,\ldots,N-1$, we conclude as above 
that $\gamma_2$ is the only left zero of the polynomial $L_{\gamma_1}f$ and that 
$\gamma_2=\gamma_2^\prime$. We subsequently get $\gamma_j=\gamma_j^\prime$ so that \eqref{5.1} is 
indeed a unique factorization of $f$. This completes the proof of $(1)\Rightarrow (4)$.

\smallskip

To prove $(4)\Rightarrow (5)$, observe that the uniqueness of \eqref{5.1} implies that any 
left divisor $h$ of $f$ is of the form $h(z)=(z-\gamma_1)(z-\gamma_2)\cdots(z-\gamma_n)$ for some 
$n\le N$ (to see this, it suffices to compare factorization $f=hg$ with \eqref{5.1}). Therefore,
for any two proper left divisors $h=\bp_{\gamma_1}\cdots\bp_{\gamma_n}$ and 
$\widetilde{h}=\bp_{\gamma_1}\cdots\bp_{\gamma_k}$ ($n\le k<N$) of $f$, we have
$\langle h\rangle_{\bf r}\cap \langle \widetilde h\rangle_{\bf r}=\langle \widetilde h\rangle_{\bf r}
\neq \langle f\rangle_{\bf r}$. Therefore, the ideal $\langle f\rangle_{\bf r}$	is irreducible.

\smallskip

To prove $(5)\Rightarrow (2)$, let us assume that the ideal $\langle f\rangle_{\bf r}$ 
irreducible. If $f$ has zeros in more than one conjugacy class, then $f$ is equal to the 
{\bf lrcm} of its left spherical divisors, by Theorem \ref{T:1.1}. 
Since in this case, each left spherical divisor of $f$ is a proper left divisor of $f$, it follows 
that the ideal $\langle f\rangle_{\bf r}$ is not irreducible which is a contradiction. Therefore 
$\mathcal Z(f)\subset [\alpha]$ for some $\alpha\in\bH\backslash\R$. If $f=\cX_{[\alpha]}^\kappa$ for $\kappa\ge 
1$, then the 
polynomials $g=\bp_{\alpha}^\kappa$ and $h=\bp_{\overline{\alpha}}^\kappa$ are proper left divisors 
of $f$ and their {\bf lrcm} equals $f$. Therefore, $\langle f\rangle_{\bf r}$ is 
not irreducible which contradicts the current assumption. It now follows from Corollary \ref{C:4.5}
that $f$ is necessarily of the form \eqref{4.14}. If $\kappa=m_s([\alpha];f)>0$, then  
the polynomials
$$
g=\bp_{\alpha_1}\cdots\bp_{\alpha_{n-1}}\bp_{\alpha_{n}}^{\kappa+1}\quad\mbox{and}\quad
h=\bp_{\overline{\alpha}_1}^\kappa
$$
are proper left divisors of $f$ and their least right common multiple equals $f$
(the details are furnished in Lemma \ref{L:5.4} below; see also Remark \ref{R:5.5}).  We again 
conclude that $\langle 
f\rangle_{\bf r}$ is not irreducible which contradicts the current 
assumption. Therefore, $\kappa=0$ in representation \eqref{4.14} and $f$ has a unique left zero.

\smallskip

We have verified implications $(1)\Leftrightarrow(2)\Rightarrow (4)\Rightarrow (5)\Rightarrow (2)$.
Implications $(1)\Leftrightarrow(3)\Rightarrow (4)\Rightarrow (6)\Rightarrow (3)$ are verified 
in much the same way.\qed

\medskip
\noindent
{\bf Notation:} In what follows, we will write ${\mathcal P}_V:=\left\{f\in\bH[z]: \; \cZ(f)\subset V\right\}$
for the set of polynomials having all zeros in $V$, and we will denote by $\mathcal{IP}_V$ the set of  
(indecomposable) polynomials having one left and one right zero in $V$.
\begin{lemma}
Let $g,h\in\mathcal{IP}_V$ ($\deg (g)=n\ge k= \deg (h)$) be given in the form
\begin{equation}
g=\bp_{\alpha_1}\cdots\bp_{\alpha_n},\quad h=\bp_{\beta_1}\cdots\bp_{\beta_k}
\label{5.2}   
\end{equation}
where $\balpha=(\alpha_1,\ldots,\alpha_n)$ and $\bbeta=(\beta_1,\ldots,\beta_k)$
are two spherical chains from the conjugacy class $V$.
If $g$ and $h$ are left coprime (i.e., if $\alpha_1\neq \beta_1$), then 
\begin{equation}
f:={\bf lrcm}(g,h)=\left\{\begin{array}{ccc} \cX_V^k, &\mbox{if} \; \;  & n=k,\\
\cX_V^k\bp_{\alpha_1}\bp_{\alpha_2}\cdots \bp_{\alpha_{n-k}}, &\mbox{if} \; \;  & n>k.
\end{array}\right.
\label{5.3}
\end{equation}
Similarly, if $g$ and $h$ are right coprime (i.e., if $\alpha_n\neq \beta_k$), then 
\begin{equation}
\widetilde{f}:={\bf llcm}(g,h)=\left\{\begin{array}{ccc} \cX_V^k, &\mbox{if}\; \; & n=k,\\
\cX_V^k\bp_{\alpha_{k+1}}\bp_{\alpha_{k+2}}\cdots \bp_{\alpha_{n}}, &\mbox{if}\; \; & n>k.
\end{array}\right.
\label{5.4}
\end{equation}
\label{L:5.2}
\end{lemma}
{\bf Proof:} Since $f$ is a right common multiple of $g$ and $h$ we have from \eqref{5.2}
\begin{equation}
f=\bp_{\alpha_1}\cdots\bp_{\alpha_{n-1}}
\bp_{\alpha_n}p=\bp_{\beta_1}\cdots\bp_{\beta_{k-1}}\bp_{\beta_k}q\quad\mbox{for some}\quad 
p,q\in\bH[z].
\label{5.5}
\end{equation}
We then use \eqref{2.17} to evaluate the latter representations 
at $\alpha_1$ and at $\beta_1$ on the left:
\begin{align*}
f^{\bl}(\beta_1)&=(\beta_1-\alpha_1)(\overline{\alpha}_1-\alpha_2)\cdots 
(\overline{\alpha}_{n-1}-\alpha_{n}) p^{\bl}(\overline{\alpha}_n)=0,\\
f^{\bl}(\alpha_1)&=0=(\alpha_1-\beta_1)(\overline{\beta}_1-\beta_2)\cdots 
(\overline{\beta}_{k-1}-\beta_{k})q^{\bl}(\overline{\beta}_k).
\end{align*}
Since $\alpha_1\neq \beta_1$, $\alpha_{j+1}\neq \overline{\alpha}_j$ and $\beta_{j+1}\neq 
\overline{\beta}_j$, the latter equalities imply 
$p^{\bl}(\overline{\alpha}_n)=q^{\bl}(\overline{\beta}_k)=0$. By \eqref{2.6},
$p=\bp_{\overline{\alpha}_n}p_1$ and $q=\bp_{\overline{\beta}_k}q_1$
for some $p_1,q_1\in\bH[z]$.
Substituting the latter factorizations into \eqref{5.5} gives
\begin{align*}
f&=\bp_{\alpha_1}\cdots\bp_{\alpha_{n-1}}\bp_{\alpha_n}\bp_{\overline{\alpha}_n}p_1=
\bp_{\alpha_1}\cdots\bp_{\alpha_{n-1}}p_1\cX_V\\
&=\bp_{\beta_1}\cdots\bp_{\beta_{k-1}}\bp_{\beta_k}\bp_{\overline{\beta}_k}q_1=
\bp_{\beta_1}\cdots\bp_{\beta_{k-1}}q_1\cX_V.
\end{align*}
Applying the spherical backward shift $S_V$ \eqref{4.9b} to the latter equalities gives
$$
S_{V}f=\bp_{\alpha_1}\cdots\bp_{\alpha_{n-1}}p_1=\bp_{\beta_1}\cdots\bp_{\beta_{k-1}}q_1.
$$
Repeating the preceding argument $k-1$ more times we get polynomials $p_1,\ldots,p_k$ and 
$q_1,\ldots,q_k$ such that $p_{j}=\bp_{\alpha_{n-j}}p_{j+1}$ and $q_{j}=\bp_{\alpha_{k-j}}q_{j+1}$, 
and eventual equalities
\begin{equation}
S_V^kf=p_k=q_k\; \; (\mbox{if} \; \; n=k)\quad\mbox{or}\quad 
S_V^kf=\bp_{\alpha_1}\cdots\bp_{\alpha_{n-k}}p_k=q_k \; \; (\mbox{if} \; \; n>k).
\label{5.6}
\end{equation}
For $f$ to be a right common multiple of the minimally possible degree, it is necessary and 
sufficient that $p_k\equiv 1$, and then $f$ is recovered from \eqref{5.6} as in \eqref{5.3}.
The formula \eqref{5.4} is justified quite similarly.\qed
\begin{corollary}
Let $g_1,\ldots,g_m\in\mathcal{IP}_V$ be such that 
$\deg(g_1)\ge \deg(g_2)\ge \ldots \ge \deg(g_m)$. 
If $g_1,\ldots,g_m$ are pairwise left (right) coprime, then 
${\bf lrcm}(g_1,g_2,\ldots,g_m)={\bf lrcm}(g_1,g_2)$ (respectively, 
${\bf llcm}(g_1,g_2,\ldots,g_m)={\bf llcm}(g_1,g_2)$).
\label{C:5.3}
\end{corollary}
{\bf Proof:} Let $\deg(g_j)=d_j$. By Lemma \ref{L:5.2}, $f:={\bf lrcm}(g_1,g_2)=\cX_V^{d_2}h$ 
for some $h\in \mathcal{IP}_V$ with $\deg(h)=d_1-d_2$.
Since $g_jg_j^\sharp=\cX_V^{d_j}$, it follows that $\cX_V^{d_2}$ is a right common multiple of 
$g_2,\ldots,g_m$. Thus, $f$ is a  right common multiple of $g_2,\ldots,g_m$ and the 
least right common multiple of $g_1$ and $g_2$. Therefore, $f={\bf lrcm}(g_1,g_2,\ldots,g_m)$,
which proves the first statement. The second statement is verified similarly.\qed

\medskip

By Corollary \ref{C:4.5}, any monic polynomial $f\in\mathcal P_V$
can be uniquely represented as the product $f=\cX_V^k p$ 
of a polynomial $p\in\mathcal {IP}_V$ and a power of 
the characteristic polynomial $\cX_V$. Formula \eqref{5.3} (along with Theorem \ref{T:5.1}) 
tells us that the ${\bf lrcm}(g,h)$ of two relatively prime polynomials $g,h\in\mathcal{IP}_V$ 
($\deg (g)\ge \deg (h)$) is of exactly the same form, where $k=\deg (h)$ and $p$ is the (unique) 
left divisor of $g$ of degree 
$$
\deg (p)=\deg(g)-\deg(h)=\deg (f)-2\deg(h).
$$ 
We thus arrive at the following result.
\begin{lemma}
Any polynomial $f\in\mathcal P_V$ with $m_s(V;f)=k\ge 1$ can be represented as 
\begin{equation}
f={\bf lrcm}(g,h),\qquad g,h\in\mathcal {IP}_V, \; \; \cZ_{\boldsymbol\ell}(g)\neq 
\cZ_{\boldsymbol\ell}(h),
\label{5.7}   
\end{equation}
and for any such representation, $\deg(h)=k$ and $\deg(g)=\deg(f)-k$. Moreover, 
\begin{enumerate}
\item If $f=\cX_V^k$, then \eqref{5.7} holds for any $g,h\in\mathcal {IP}_V$ 
with $\deg(h)=\deg(g)=k$.
\item If $f=\cX_V^k p$ with $p=\bp_{\alpha_1}\ldots \bp_{\alpha_n}\in\mathcal {IP}_V$, then
all pairs $(g,h)$ giving rise to representation \eqref{5.7} are characterized by the properties
\begin{enumerate}
\item $\deg(h)=k$ and $\cZ_{\boldsymbol\ell}(h)\neq \{\alpha_1\}$;
\item $g=pq\quad\mbox{for some}\quad q\in\mathcal {IP}_V$ such that $\deg (q)=k$ and 
$\cZ_{\boldsymbol\ell}(g)\neq \{\overline{\alpha}_n\}$.
\end{enumerate}
\end{enumerate}
\label{L:5.4}
\end{lemma}
\begin{Remark}
In Lemma \ref{L:5.4}, one can choose $h=\bp_{\overline{\alpha}}^k$ and $g=\bp_{\alpha}^k$ 
for any fixed $\alpha\in V$ in case (1)  and $h=\bp_{\overline{\alpha}_1}^k$ and 
$q=\bp_{\alpha_n}^k$ in case (2). Observe that these particular choices were used 
in the proof of Theorem \ref{T:5.1}.
\label{R:5.5}
\end{Remark}
We now formulate a more detailed version of Theorem \ref{T:1.2}. Recall that noncommutative 
polynomials $p_1,\ldots,p_n$ are said to be {\em left (right) relatively prime}
if each polynomial $p_k$ has no common left (right) zeros with the least right (left) common multiple of all other
polynomials. We also recall
that the spherical divisors $D^f_{{\boldsymbol \ell}, V_j}$ and $D^f_{{\bf r}, V_j}$ of a polynomial $f$ are unique
by Theorem \ref{1.1}.
\begin{theorem}
Let $f\in\bH[z]$ be a monic polynomial with spherical zeros $V_1,\ldots,V_m$ 
and isolated zeros contained in conjugacy classes $V_{m+1},\ldots,V_n$.
There exist two sets $\Pi=\{p_i\}_{i=1}^{m+n}$ and 
$\widetilde{\Pi}:=\{\widetilde{p}_i\}_{i=1}^{m+n}$ 
of relatively prime indecomposable polynomials such that 
\begin{equation}
\langle f\rangle_{\bf r} =\bigcap_{i=1}^{m+n}\langle p_i\rangle_{\bf r}\quad\mbox{and}\quad
\langle f\rangle_{\boldsymbol\ell}=\bigcap_{i=1}^{m+n}\langle
\widetilde{p}_i\rangle_{\boldsymbol\ell}.
\label{5.8}
\end{equation}
Representations \eqref{5.8} are unique in the following sense:
\begin{enumerate}
\item For each  $j=m+1,\ldots,n$, the set $\Pi$ (resp., $\widetilde{\Pi}$) contains exactly one 
polynomial  from $\mathcal {IP}_{V_j}$
(for each $j=m+1,\ldots,n$) which is equal to $D^f_{{\boldsymbol \ell}, V_j}$ 
(resp., $D^f_{{\bf r}, V_j}$). 
\item For each  $j=1,\ldots,m$, the set $\Pi$ (resp., $\widetilde{\Pi}$) contains exactly 
two polynomials from $\mathcal {IP}_{V_j}$, and the least right (left) common multiple 
of these polynomials is equal to $D^f_{{\boldsymbol \ell}, V_j}$
(resp., $D^f_{{\bf r}, V_j}$).
\end{enumerate}
\label{T:5.6}
\end{theorem}
{\bf Proof:} The existence of representations \eqref{5.8} follows from 
Theorem \ref{T:1.1} and Lemma \ref{L:5.4}. Let us assume that
\begin{equation}
\langle f\rangle_{\bf r} =\bigcap_{i=1}^{M}\langle p_i\rangle_{\bf r}
\label{5.9}
\end{equation}
for a relatively prime collection $\Pi=\{p_i: \; 1\le i\le M\}$ of indecomposable polynomials.
In particular, each polynomial $p_i\in\Pi$ is relatively prime with the {\bf lrcm}
of all other polynomials in $\Pi$. By Corollary \ref{C:5.3}, it follows that 
$\Pi$ contains at most two polynomials with zeros in the same conjugacy class. 
Moreover, if $\Pi$ contains two polynomials in $\mathcal{IP}_{V_j}$, then 
$V_j$ is a spherical zero of $f$, i.e., $j\in\{1,\ldots,m\}$. Otherwise, $V_j$ contains 
isolated zeros of $f$, i.e., $j\in\{m+1,\ldots,n\}$. Therefore, $M=2m+(n-m)=m+n$.
Let us define $F_j\in\mathcal {P}_{V_j}$ as the {\bf lrcm} of two elements in $\Pi\cap \mathcal 
{IP}_{V_j}$ for $j=1,\ldots,m$ or as a unique element in $\Pi\cap \mathcal {IP}_{V_j}$
for $j\ge m$. Then we have from \eqref{5.9},
$$
\langle f\rangle_{\bf r} =\bigcap_{j=1}^{n}\langle F_j\rangle_{\bf r}, \qquad F_j\in\mathcal 
{P}_{V_j}.
$$
As we know from Theorem \ref{T:1.1}, the latter representation implies 
$F_j=D^f_{{\boldsymbol \ell}, V_j}$ for $j=1,\ldots,n$. This completes the proof 
of the part concerning the first representation in \eqref{5.8}. The dual part is
verified in much the same way.\qed

\section{Least common multiples}
\setcounter{equation}{0}

In the two previous sections, we represented a given polynomial $f\in\bH[z]$
as the least common multiple of its divisors of certain type. Now
we address the converse problem. 
\begin{problem}
Given a finite collection $\Pi=\{g_1,\ldots,g_m\}\subset\bH[z]$, construct explicitly
${\bf lrcm}(g_1,\ldots,g_m)$ and ${\bf llcm}(g_1,\ldots,g_m)$. 
\label{P:6.1}
\end{problem}
Lemma \ref{L:5.2} and Corollary \ref{C:5.3} settled the case where $\Pi\subset\mathcal 
{IP}_V$ consists of pairwise 
coprime polynomials. The next result removes the coprimeness assumption. 
Throughout the section, we will be dealing only with right common multiples
and consequently, with the left zero structure. The dual statements are analogous
and will be omitted.
\begin{lemma}
Given $g_1,\ldots,g_m\in\mathcal{IP}_V$ such that  
$\deg(g_1)\ge \deg(g_2)\ge \ldots \ge \deg(g_m)$, let 
$g_1={\displaystyle\prod_{i=1}^{\substack{\curvearrowright \\ n}}\bp_{\alpha_i}}:=
\bp_{\alpha_1}\bp_{\alpha_2}\cdots \bp_{\alpha_n}$ and let
\begin{equation}
g_j=p_jh_j,\quad\mbox{where}\quad p_j={\bf glcd}(g_j,g_1)\quad\mbox{for}\quad
j=2,\ldots,m.
\label{6.1}
\end{equation}
Then 
\begin{equation}
{\bf lrcm}(g_1,g_2,\ldots,g_m)=\cX_V^k \bp_{\alpha_1}\bp_{\alpha_2}\cdots 
\bp_{\alpha_{n-k}},\quad\mbox{where}\quad k=\max_{2\le j\le m} \deg (h_j).
\label{6.2}
\end{equation}
\label{L:6.2}
\end{lemma}
{\bf Proof:} Since $g_1$ is indecomposable, its left divisor $p_j$ is of the 
form $p_j={\displaystyle\prod_{i=1}^{\substack{\curvearrowright \\ 
\deg(p_j)}}\bp_{\alpha_i}}$,
by property (4) in 
Theorem \ref{T:5.1}. Therefore, for a fixed $j\in\{2,\ldots,m\}$, we have 
\begin{align}
{\bf lrcm}(g_j,g_1)&=p_j\cdot {\bf lrcm}\bigg(h_j,
\prod_{i=\deg(p_j)+1}^{\substack{\curvearrowright \\ n}}\bp_{\alpha_i}\bigg)\notag\\
&=\cX_V^{\deg(h_j)}\cdot p_j\cdot \prod_{i=\deg(p_j)+1}^{\substack{\curvearrowright \\ 
n-\deg(h_j)}}\bp_{\alpha_i}\notag\\
&=\cX_V^{\deg(h_j)}\cdot 
\prod_{i=1}^{\substack{\curvearrowright \\
n-\deg(h_j)}}\bp_{\alpha_i}={\bf 
lrcm}\left(\bp_{\overline{\alpha}_1}^{\deg(h_j)},g_1\right).\label{6.3}
\end{align}
The first equality in the latter calculation follows from \eqref{6.1}, the second 
follows by applying Lemma \ref{L:5.2} to left coprime polynomials 
$h_j$ and $\bp_{\alpha_{\deg(p_j)}}\cdots \bp_{\alpha_n}$ and since the polynomial $\cX_V$ 
is  real, the third equality follows from factorization of $p_j$, and the last equality 
follows from Lemma \ref{L:5.4} (see also Remark \ref{R:5.5}). 
We now get \eqref{6.2}:
\begin{align*}
{\bf lrcm}(g_1,g_2,\ldots,g_m)&={\bf lrcm}(g_1, {\bf lrcm}(g_2,g_1),\ldots,
{\bf lrcm}(g_m,g_1))\\
&={\bf lrcm}\left(g_1, {\bf lrcm}\left(\bp_{\overline{\alpha}_1}^{\deg(h_2)},g_1\right),
\ldots,{\bf lrcm}\left(\bp_{\overline{\alpha}_1}^{\deg(h_m)},g_1\right)\right)\\
&={\bf lrcm}\left(g_1, \bp_{\overline{\alpha}_1}^{\deg(h_2)},\ldots,
\bp_{\overline{\alpha}_1}^{\deg(h_m)}\right)\\
&={\bf lrcm}\left(g_1, \bp_{\overline{\alpha}_1}^{k}\right)=
\cX_V^k \bp_{\alpha_1}\bp_{\alpha_2}\cdots \bp_{\alpha_{n-k}},
\end{align*}
where the first and the third equalities are self-evident, the second equality holds 
due to \eqref{6.3}, the fourth equality holds due to the choice \eqref{6.2} of $k$, and 
the last equality holds by Lemma \ref{L:5.2} applied to left coprime polynomials 
$g_1$ and $\bp_{\overline{\alpha}_1}^{k}$.\qed

\medskip

Observe that if $g_2,\ldots,g_m$ are left coprime with $g_1$, then $p_j\equiv 1$ for
$j=2,\ldots,m$ and Corollary \ref{C:5.3} follows from Lemma \ref{L:6.5}.
We next consider Problem \ref{P:6.1} for indecomposable polynomials 
having zeros in distinct conjugacy classes. The case where all polynomials are linear
($\Pi=\{\bp_{\gamma_1},\ldots,\bp_{\gamma_n}\}$), 
has been known for a while; see e.g., \cite{cm}. 
To handle the general case, we need the following preliminary result.
\begin{lemma}
Let $V$ be a conjugacy class, let $F=\bp_{\alpha_1}\cdots \bp_{\alpha_k}\in\mathcal {IP}_{V}$
and let $Q\in\bH[z]$ be such that $\cZ(Q)\cap V=\emptyset$. Define
$\widetilde{\alpha}_1,\ldots,\widetilde{\alpha}_k\in\bH$
and $Q_0, Q_{1},\ldots,Q_k\in\bH[z]$ by
\begin{equation}
Q_0=Q,\quad \widetilde{\alpha}_j=Q_{j-1}^{\bl}(\alpha_j)^{-1}\alpha_j Q_{j-1}^{\bl}(\alpha_j),  
\quad Q_{j}=L_{\alpha_j}\left(Q_{j-1}\bp_{\widetilde{\alpha}_j}\right)
\label{6.6}   
\end{equation}
for $j=1,\ldots,k$. Then $\cZ(Q_k)\cap V=\emptyset$ and
\begin{equation}
{\bf lrcm}(F,Q)=FQ_k=
Q\bp_{\widetilde{\alpha}_1}\bp_{\widetilde{\alpha}_2}\cdots
\bp_{\widetilde{\alpha}_k}.
\label{6.7}  
\end{equation}
\label{L:6.5}
\end{lemma}
{\bf Proof:} Since $\cZ(F)\subset V$ and $\cZ(Q)\cap V=\emptyset$, it follows
(see \cite[Proposition 4.2]{lamler} for the proof) that
$\deg({\bf lrcm}(F,Q))=\deg(F)+\deg(Q)=\deg(F)+k$. Therefore, it suffices to find 
a right common multiple of polynomials $F$ and $Q$ of degree equal $\deg(F)+k$.

\smallskip

Based on the second statement in Lemma \ref{L:4.7} and definitions
\eqref{6.6}, an induction argument shows that
$$
\bp_{\alpha_j}Q_{j}=Q_{j-1}\bp_{\widetilde\alpha_{j}}\quad\mbox{and}\quad \cZ(Q_{j})\cap V=\emptyset
\quad\mbox{for}\quad j=1,\ldots,k.
$$
Therefore, we get recursively
\begin{align*}
FQ_k&=\bp_{\alpha_1}\bp_{\alpha_2}\cdots \bp_{\alpha_{k-1}}\bp_{\alpha_{k}}Q_k\\
&=\bp_{\alpha_1}\bp_{\alpha_2}\cdots\bp_{\alpha_{k-1}}Q_{k-1}\bp_{\widetilde{\alpha}_k}\\
&=\bp_{\alpha_1}\cdots\bp_{\alpha_{k-2}}Q_{k-2}\bp_{\widetilde{\alpha}_{k-1}}
\bp_{\widetilde{\alpha}_{k}}=\ldots =Q_{0}\bp_{\widetilde{\alpha}_1}\bp_{\widetilde{\alpha}_{2}}
\cdots\bp_{\widetilde{\alpha}_{k}}=Q\bp_{\widetilde{\alpha}_1}\bp_{\widetilde{\alpha}_{2}}\cdots 
\bp_{\widetilde{\alpha}_{k}},
 \end{align*}
from which we conclude that the polynomial $FQ_k$ is a right common multiple of $F$ and $Q$. 
It is clear that  its degree equals  $\deg(F)+k$, so that \eqref{6.7} follows.\qed

\medskip

The next algorithm produces the {\bf lrcm} of $n$ indecomposable polynomials 
\begin{equation}
F_i=\bp_{\alpha_{i,1}}\bp_{\alpha_{i,2}}
\cdots \bp_{\alpha_{i,k_i}}\in\mathcal{IP}_{V_i} \quad (i=1,\ldots,n)
\label{6.8}      
\end{equation}   
with zeros in $n$ distinct conjugacy classes $V_1,\ldots,V_n\subset\bH$.
\begin{algorithm} Given polynomials \eqref{6.8},
\begin{enumerate}
\item Let $P_1:=F_1=\bp_{\alpha_{1,1}}\bp_{\alpha_{1,2}}
\cdots \bp_{\alpha_{1,k_1}}$ and let $i:=2$.
\item Let $Q_0:=P_{i}$ and  perform the recursion
\begin{equation}
\widetilde{\alpha}_{i+1,j}=Q_{j-1}^{\bl}(\alpha_{i+1,j})^{-1}\alpha_{i+1,j} 
Q_{j-1}^{\bl}(\alpha_{i+1,j}),
\quad Q_{j}=L_{\alpha_{i+1,j}}\left(Q_{j-1}\bp_{\widetilde{\alpha}_{i+1,j}}\right)
\label{6.9}
\end{equation}
for $j=1,\ldots,k_i$.
\item Let $P_{i+1}:=P_{i}Q_{k_i}$. If $i<n-1$, then let $i:=i+1$ and go to {\rm (2)}.
If $i=n-1$, proceed to {\rm (4)}
\item The polynomial $P_n$ is equal to the ${\bf lrcm}(F_1,\ldots,F_n)$.
\end{enumerate}
\label{A:6.6}  
\end{algorithm}
To justify the algorithm, let us assume that $P_i={\bf lrcm}(F_1,\ldots,F_{i})$.
Then it follows from Lemma \ref{L:6.5} that 
\begin{align*}
P_{i+1}:=P_{i}Q_{k_i}&={\bf lrcm}(P_{i},F_{i+1})\\
&={\bf lrcm}({\bf lrcm}(F_1,\ldots,F_{i}), \, F_{i+1})={\bf lrcm}(F_1,\ldots,F_{i+1}).
\end{align*}
Since $P_1=F_1={\bf lrcm}(F_1)$, we conclude by induction that $P_i={\bf lrcm}(F_1,\ldots,F_{i})$
for all $i=1,\ldots,n$.\qed
\begin{Example} We illustrate Algorithm \ref{A:6.6} by constructing ${\bf
lrcm}(\bp_\alpha^2,\bp_\beta^2)$ for two non-real quaternions $\alpha\not\sim\beta$.
We let $P_1=\bp_\alpha^2$, and perform two steps of recursion \eqref{6.9}:
\begin{align}
Q_0&=P_1=\bp_\alpha^2,\quad Q_0^{\bl}(\beta)=\beta^2-2\beta\alpha+\alpha^2,\notag\\
\beta_1&=(\beta^2-2\beta\alpha+\alpha^2)^{-1}\beta (\beta^2-2\beta\alpha+\alpha^2),\label{6.10}\\
Q_1&=L_{\beta}\left(\bp_\alpha^2 \bp_{\beta_1}\right)=z^2+(\beta-2\alpha-\beta_1)z+
\beta^2-2\beta\alpha+2\alpha\beta_1-\beta\beta_1+\alpha^2,\notag\\
Q_1^{\bl}(\beta)&=3\beta^2-4\beta\alpha+\alpha^2+2(\alpha-\beta)\beta_1,\notag\\
\beta_2&=(3\beta^2-4\beta\alpha+\alpha^2+2(\alpha-\beta)\beta_1)^{-1}
\beta(3\beta^2-4\beta\alpha+\alpha^2+2(\alpha-\beta)\beta_1).\notag
\end{align}
Now we can write the answer:
\begin{equation}
G(z):={\bf lrcm}(\bp_\alpha^2,\bp_\beta^2)=(z-\alpha)^2(z-\beta_1)(z-\beta_2).
\label{6.11}
\end{equation}
If we choose in the previous example, $\alpha={\bf i}$ and $\beta=1+{\bf j}$, then formulas
\eqref{6.10} give
\begin{align*} 
\beta_1&=(-1-2{\bf i}+2{\bf j}+2{\bf k})^{-1}(1+{\bf j})(-1-2{\bf i}+2{\bf j}+2{\bf k})=
1-\frac{12{\bf i}+3{\bf j}-4{\bf k}}{13},\\
\beta_2&=(-21-10{\bf i}+50{\bf j}+14{\bf k})^{-1}(1+{\bf j})(-21-10{\bf i}+50{\bf j}+14{\bf k})\\
&=1+\frac{-1588{\bf i}+2645{\bf j}+980{\bf k}}{3237},
\end{align*}
and we get the {\bf lrcm} of polynomials $(z-{\bf i})^2$ and $(z-1-{\bf j})^2$ using formula
\eqref{6.11}.
\label{E:6.7}  
\end{Example}
We would like to stress that to ensure  $m_{\boldsymbol\ell}(\beta;G)=2$ for $G$ of the form 
\eqref{6.11}, we cannot choose $\beta_2=\beta_1$;  actually, it can be shown that ${\bf 
lrcm}(\bp_\alpha^2,\bp_\beta^2)=\bp_\alpha^2 \bp_{\beta_1}^2$ (with $\beta_1$ as in 
\eqref{6.10}) if and only if $\alpha\beta=\beta\alpha$ in which case we also have 
$\beta_1=\beta$.  Thus, finding least common multiples of special irreducible polynomials 
$\bp_{\alpha_i}^{k_i}$ is  the computational problem of about the same complexity as the 
generic one. 

\medskip

We now present several algorithms based on Algorithm \ref{A:6.6}. 
The first algorithm constructs a polynomial with prescribed left spherical divisors.
The input is the collection of polynomials $F_j\in\mathcal {P}_{V_j}$
($j=1,\ldots,m$), and the output is $f={\bf lrcm}(F_1,\ldots,F_n)$, which, according 
to Theorem \ref{T:1.1}, is a unique monic polynomial such that 
\begin{equation}
\cZ(f)\subset\bigcup_{j=1}^m V_j\quad\mbox{and}\quad
D^f_{{\boldsymbol \ell}, V_j}=F_j\quad\mbox{for}\quad j=1,\ldots,m.
\label{8.1a}
\end{equation}
By Corollary \ref{C:4.5}, each $F_j$ is 
either of the form $F_j=\bp_{x_j}^{k_j}$ (if $V_j=\{x_j\}\subset \R$) or, otherwise,
$F_j=\cX_{V_j}^{\kappa_j}p_j$, where either $p_j\equiv 1$ or $p_{j}\in\mathcal {IP}_{V_j}$. 
\begin{algorithm}
Let $F_1,\ldots,F_m$ ($F_j\in\mathcal{P}_{V_j}$) be arranged so that
$$
F_j=\cX_{V_j}^{\kappa_j}p_j \; (1\le j\le m_1)\quad F_j=\cX_{V_j}^{\kappa_j} \;  (m_1<j\le m_2);\quad
F_j=\bp_{x_j}^{k_j} \; (m_2<j\le m),
$$
where $p_j\in\mathcal{IP}_{V_j}$ for $j=1,\ldots,m_1$. 
Use Algorithm \ref{A:6.6} to compute $G={\bf lrcm}(p_{1},\ldots,p_{m_1})$ and write
\begin{equation}
{\bf lrcm}(F_1,\ldots,F_n)=G\cdot \bigg(\prod_{j=1}^{m_2}\cX_{V_j}^{\kappa_j}\bigg)\cdot 
\bigg(\prod_{j=m_1+1}^{m_2}\bp_{x_j}^{\kappa_j}\bigg).
\label{8.2}
\end{equation}
\label{A:8.1}
\end{algorithm}
{\bf Proof:} Algorithm \ref{A:6.6} applies since $p_1,\ldots,p_{m_1}$ are indecomposable polynomials with zeros in 
$m_1$ distinct conjugacy classes. The rest is clear since $\bp_{x_j}^{k_j}$ and $\cX_{V_j}^{\kappa_j}$ are real 
polynomials and since 
${\bf lrcm}(F_1,\ldots,F_{m_1})=\left(\prod_{j=1}^{m_1}\cX_{V_j}^{\kappa_j}\right)\cdot {\bf 
lrcm}(p_1,\ldots,p_{m_1})$.\qed

\medskip

We next consider the case where all polynomials have zeros in the same (non-real) conjugacy class $V\subset\bH$. 
Without loss of generality we may assume that at most one of given polynomials is real (i.e., of the form
$\cX_V^\kappa$) -- if there are several, then removing all of them but the one of highest degree will 
not
affect the {\bf lrcm}.
\begin{algorithm}
Given polynomials $F_1=\cX_V^{\kappa_1}$ and $F_j=\cX_V^{\kappa_j} p_j$ ($j=2,\ldots,m$),
where $p_j=\bp_{\alpha_{j,1}}\bp_{\alpha_{j,2}}\cdots \bp_{\alpha_{j,k_j}}\in\mathcal{IP}_V$,
\begin{enumerate}
\item Construct the polynomials $g_{j1}$ and $g_{j2}$ as follows: 
\begin{enumerate}
\item Let $g_{11}=\bp_{\alpha}^{\kappa_1}$ and $g_{12}=\bp_{\overline{\alpha}}^{\kappa_1}$ 
for any fixed $\alpha\in V$.
\item For $j=2,\ldots,m$, let $g_{j1}=p_j\bp_{\alpha_{j,k_j}}^{\kappa_j}$ and 
$g_{j2}=\bp_{\overline{\alpha}}^{\kappa_j}$.
\end{enumerate}
\item Find $f={\bf lrcm}(g_{j1},g_{j2}: \; 1\le j\le m)$ as in Lemma \ref{L:6.2}. Then we 
also have $f={\bf lrcm}(F_1,\ldots,F_{m})$.
\end{enumerate}
\label{A:8.2}
\end{algorithm}
The statement is obvious since by By Remark \ref{R:5.5}, $F_j={\bf lrcm}(g_{j1}, \, g_{j2})$ for $j=1,\ldots,m$.
Note that in Step 2, we may dismiss all polynomials $F_j$ of degree $\deg F_j\le \kappa_1$ (since 
$F_1$ is a multiple of each such polynomial). Also, if the polynomials $p_j\in\mathcal{IP}_V$ are pairwise left coprime, 
then we may choose $\alpha$ in Step 1(a) so that the polynomials $g_{ji}$ are all pairwise coprime in which case 
Step 2 simplifies to 

\smallskip

$2^\prime.$ {\em Pick (any) two polynomials in the set $\{g_{j1},g_{j2}: \; 1\le j\le m\}$ with highest degrees
and find their least right common multiple using formula \eqref{5.3}.}

\medskip

We finally present the algorithm that produces the {\bf lrcm} of any finite collection of 
quaternion polynomials. The construction is based on the fact that the left spherical divisor of the 
${\bf lrcm}(g_1,\ldots,g_m)$ associated with a conjugacy class $V$ is equal to the 
{\bf lrcm} of left spherical divisors of $g_1,\ldots,g_m$ associated with $V$.
\begin{algorithm}
Given polynomials $g_1,\ldots,g_m\in\bH[z]$,
\begin{enumerate}
\item Find all conjugacy classes $V_1,\ldots,V_n$ containing at least one left zero of at least one polynomial from 
the set.
\item Use Algorithm \ref{A:4.3} to find all spherical divisors $D_{\boldsymbol\ell,V_i}^{g_j}$ of $g_j$ for 
$j=1,\ldots,m$. If $\cZ(g_j)\cap V_i=\emptyset$, let $D_{\boldsymbol\ell,V_i}^{g_j}\equiv 1$.
\item For each $i=1,\ldots,n$, use Algorithm \ref{A:8.2} to construct 
$F_i={\bf lrcm}(D_{\boldsymbol\ell,V_i}^{g_1}, 
\ldots,D_{\boldsymbol\ell,V_i}^{g_m})$.
\item Use Algorithm \ref{A:8.1} to construct $f={\bf lrcm}(F_1,\ldots,F_n)$.
\end{enumerate}
Then we also have $f={\bf lrcm}(g_1,\ldots,g_m)$.
\label{A:8.3} 
\end{algorithm}

\section{Formal power series over quaternions}
\setcounter{equation}{0}
As in the commutative case, certain results concerning quaternion polynomials can be extended
to formal power series over $\bH$ (see e.g., \cite{gs}, \cite{gss}, \cite{abcs}). We are particularly interested 
in power series for which left and right evaluation functionals make sense.
We denote by $\mathbb B=\left\{\alpha\in\bH: \,  |\alpha|<1 \right\}$ the open unit ball in $\bH$,  
and introduce the ring
$$
\mathcal H=\bigg\{f(z)=\sum_{j=0}^\infty f_jz^j: \; \limsup \sqrt[k]{|f_k|}\le 1\bigg\}.
$$
Observe that for any $f(z)=\sum f_kz^k$ in $\mathcal H_R$ and any $\alpha\in\mathbb B$, the
quaternion series $\sum_{k=0}^\infty\alpha^k f_k$ and $\sum_{k=0}^\infty f_k\alpha^k$
converge absolutely, so the evaluation formulas \eqref{2.2} (with $m=\infty$) and therefore, 
the notions of left and right zeros (within $\mathbb B$) make sense. Furthermore, the power series $L_\alpha f$ and 
$R_\alpha f$ (defined as in \eqref{2.3} but with $m=\infty$) are also in 
$\mathcal H$. Equalities \eqref{2.1} and therefore, equivalences \eqref{2.5} hold true in $\mathcal H$ as well as 
evaluation formulas \eqref{2.17}. The conjugate power series $f^\sharp$ is defined as in \eqref{2.8a} (with $m=\infty$)
and relations \eqref{2.10} hold for $f,g\in\mathcal H$. Moreover, the real power series $ff^\sharp\in\mathbb R[[z]]$
belongs to $\mathcal H$ and (by the complex uniqueness theorem) has countably many spherical zeroes in $\mathbb B$ 
of finite multiplicity each. Theorem \ref{T:2.4} extends to $\mathcal H$ as 
follows: {\em if $f\in\mathcal H$, then each conjugacy class $V\subset\cZ(ff^\sharp)$ 
either contains exactly one left and one right zero of $f$ or $V\subset\cZ_{\boldsymbol\ell}(f)\bigcap 
\cZ_{\bf r}(f)$.}
\begin{Remark}
Theorem \ref{T:2.1} holds for any $f\in\mathcal H$. 
\label{R:8.1}
\end{Remark}
Indeed, formulas \eqref{3.1} rely on representation formulas \eqref{3.3} and \eqref{3.4}, which 
hold true for monomials $f_jz^j$ and therefore, for elements in $\mathcal H$, by linearity.
Furthermore, since $m_s([\alpha],ff^\sharp)$ is finite for every $\alpha\in\mathbb B$, Algorithms \ref{A:4.1}
and \ref{A:4.3} apply to $f\in\mathcal H$; the only modification is that $Q_j$ and $\widetilde{Q}_j$ obtained via 
recursions \eqref{4.3}-\eqref{4.7} are power series from $\mathcal H$ rather than polynomials. Algorithm \ref{A:4.1}
recovers in a constructive way the following result from \cite{gs}.
\begin{proposition}
Given $f\in\mathcal H$, let $V$ be a conjugacy class such that $m_s(V,ff^\sharp)=k$.
There exist (unique) integer $\kappa\ge 0$, spherical chains $\balpha=(\alpha_1,\ldots,\alpha_{k-2\kappa})\subset V$ and 
$\widetilde\balpha=(\widetilde\alpha_1,\ldots,\widetilde\alpha_{k-2\kappa})\subset V$, and power series 
$P,\widetilde P\in\mathcal H$ having no zeros in $V$ such that 
$$
f=\cX_{V}^\kappa\bp_{\alpha_1}\bp_{\alpha_2}\cdots \bp_{\alpha_{k-2\kappa}}P=
\widetilde{P}\bp_{\widetilde{\alpha}_{k-2\kappa}}\cdots \bp_{\widetilde{\alpha}_{2}}\bp_{\widetilde{\alpha}_{1}}
\cX_{V}^\kappa.
$$
\label{P:8.1}
\end{proposition}
Although Algorithms \ref{A:4.8} and \ref{A:4.9} relating left an right zero structures 
apply to formal power series, they might be efficient only if evaluations $f^{\bl}$ and $f^{\br}$ 
admit closed formulas. 
\begin{Remark}
Proposition \ref{P:8.1} allows us to introduce
spherical divisors of a given $f\in\mathcal H$ (which are polynomials uniquely determined from $f$). Constructing
an $f\in\mathcal H$ subject to conditions \eqref{8.1a} (i.e., with finitely many prescribed spherical divisors
$D^f_{\boldsymbol\ell,V_i}$) reduces to Algorithm \ref{A:8.1} since any such $f$ is of the form $f=GH$ where $G={\bf
lrcm}(D^f_{\boldsymbol\ell,V_i},\ldots,D^f_{\boldsymbol\ell,V_m})$ and $H\in\mathcal H$ has no zeros in $\mathbb B$.
The latter construction does not seem particularly interesting unless we desire to construct an $f$ with some 
extra properties; then  the relevant question is how to choose $H$ to achieve this. In the concluding section we 
consider a problem of this sort.
\label{R:7.10}
\end{Remark}
\noindent
{\bf 7.1 Finite Blaschke products with prescribed zero structure.} We denote by
$$
{\rm H}^2=\bigg\{h(z)=\sum_{j=0}^\infty z^jh_j: \; \|h\|_{{\rm H}^2}^2:=\sum_{j=0}^\infty
|h_j|^2<\infty\bigg\}
$$
the space of elements in $\bH[[z]]$ with square summable coefficients and observe that ${\rm H}^2\subset
\mathcal H$. An element $f\in{\rm H}^2$ is called {\em bi-inner} if the equalities
\begin{equation}
\|f\cdot h\|_{{\rm H}^2}=\|h\|_{{\rm H}^2}=\|h\cdot f\|_{{\rm H}^2}
\quad\mbox{for all}\quad h\in{\rm H}^2
\label{9.3}
\end{equation}
hold.
If $f$ satisfies only the first (only the second) equality in \eqref{9.3},
it is called {\em left-inner} ({\em right-inner}). It is easily seen that 
$f$ is left-inner if and only if $f^\sharp$ is right-inner.

\smallskip

For each fixed $\alpha\in\mathbb B$, the power series
\begin{equation}
{\bf k}_{\alpha}(z)=\sum_{k=0}^\infty \alpha^kz^k\quad\mbox{belongs to ${\rm H}^2$ and $\; \|{\bf k}_{\alpha}\|_{{\rm 
H}^2}^2=\frac{1}{1-|\alpha|^2}$}.
\label{9.1}
\end{equation}
Letting $\Upsilon_{\alpha,\gamma}:=1-(\alpha+\overline{\alpha})\gamma+|\alpha|^2\gamma^2$
and observing equalities 
$$
\Upsilon_{\alpha,\gamma}\cdot \bigg(\sum_{k=0}^\infty
\gamma^k\alpha^k\bigg)=1-\gamma\overline\alpha\quad\mbox{and}\quad
\bigg(\sum_{k=0}^\infty
\alpha^k\gamma^k\bigg)\cdot\Upsilon_{\alpha,\gamma}=1-\overline\alpha\gamma\quad (|\gamma|<1),
$$
we come up with evaluation formulas 
\begin{equation}
{\bf k}_{\alpha}^{\bl}(\gamma)=\Upsilon_{\alpha,\gamma}^{-1}\cdot (1-\gamma\overline\alpha),\qquad
\quad {\bf k}_{\alpha}^{\br}(\gamma)=(1-\overline\alpha\gamma)\cdot\Upsilon_{\alpha,\gamma}^{-1}\quad \mbox{for all}\quad 
\gamma\in\mathbb B,
\label{9.2}
\end{equation}
which imply in particular, that ${\bf k}_{\alpha}$ has no zeros in $\mathbb B$. Since ${\bf 
k}_{\alpha}(z)\cdot (1-z\alpha)\equiv 1$, the power series
${\bf k}_{\alpha}$ is the formal inverse of the polynomial $1-z\alpha$. Hence,
the power series  
\begin{equation}
{\bf b}_{\alpha}(z):=\bp_{\alpha}(z)\cdot {\bf k}_{\overline\alpha}(z)=(z-\alpha)\cdot \sum_{k=0}^\infty
\overline{\alpha}^kz^k=-\alpha+(1-|\alpha|^2)z\cdot {\bf k}_{\overline\alpha}(z)
\label{9.4}
\end{equation}
can be viewed as the quaternionic analog of the Blaschke factor $\frac{z-\alpha}{1-z\overline{\alpha}}$. 
It is seen from the second representation in \eqref{9.4} that $\cZ_{\boldsymbol\ell}({\bf b}_\alpha)=
\cZ_{\bf r}({\bf b}_\alpha)=\{\alpha\}$.  In analogy to the complex case, ${\bf
b}_{\alpha}$ gives rise to two automorphisms $\gamma\mapsto {\bf b}_{\alpha}^{\bl}(\gamma)$
and $\gamma\mapsto {\bf b}_{\alpha}^{\br}(\gamma)$ of the closed unit ball $\overline{\mathbb B}$
(see \cite{heja}, \cite{bige}). The latter can be derived from the evaluation formulas
$$
{\bf b}_{\alpha}^{\bl}(\gamma)=\Upsilon_{\alpha,\gamma}^{-1}\cdot (\gamma(\alpha^2+1)-(\gamma^2+1)\alpha),
\quad {\bf k}_{\alpha}^{\br}(\gamma)=((\alpha^2+1)\gamma-\alpha(\gamma^2+1))\cdot\Upsilon_{\alpha,\gamma}^{-1}
$$
which, in turn, follow from \eqref{9.2} (with $\overline\alpha$ instead of $\alpha$) and \eqref{9.3}.
\begin{proposition} 
For each $\alpha\in\mathbb B$, the power series ${\bf b}_{\alpha}$ is bi-inner.
\label{P:8.3}
\end{proposition}
{\bf Proof:} The left-inner property ${\bf b}_{\alpha}$ has been observed in \cite{acs}. Our justification 
is similar to that  in \cite{acs} but does not rely on the Hilbert space structure of ${\rm H}^2$.
We first verify that $\|{\bf b}_{\alpha}\cdot h\|_{{\rm H}^2}=\|h\|_{{\rm H}^2}$ for every $h$ of the 
form $h(z)=d+cz^k$. Since
\begin{align*}
{\bf b}_{\alpha}(z)\cdot h(z)=&-\alpha d+(1-|\alpha|^2)\cdot \sum_{j=1}^{k-1}\overline{\alpha}^{j-1}dz^j+
((1-|\alpha|^2)\overline{\alpha}^{k-1}d-\alpha c)z^k\\
&+(1-|\alpha|^2)\cdot \sum_{j=0}^{\infty}\overline{\alpha}^{j}
(\overline{\alpha}^{k}d+c)z^{j+k+1},
\end{align*}
we have, by the definition of the ${\rm H}^2$-norm,
\begin{align}
\|{\bf b}_{\alpha}\cdot h\|^2_{{\rm H}^2}=&|\alpha|^2| d|^2+(1-|\alpha|^2)\cdot (1-|\alpha|^{2k-2})|d|^2
+|(1-|\alpha|^2)\overline{\alpha}^{k-1}d-\alpha c|^2\notag \\
&+(1-|\alpha|^2)\cdot |\overline{\alpha}^{k}d+c|^2\notag\\
=&|d|^2+|c|^2+(1-|\alpha|^2)\cdot \left(2{\rm Re}(\overline{\alpha}^{k}d\overline{c})-2{\rm Re}(
\overline{\alpha}^{k-1}d\overline{c} \, \overline{\alpha})\right).\label{9.5}
\end{align}
Applying the well-known quaternionic equality ${\rm Re}(ab)={\rm Re}(ba)$ to $a=\overline{\alpha}$ and 
$b=\overline{\alpha}^{k-1}d\overline{c}$, we conclude from \eqref{9.5} that 
$\|{\bf b}_{\alpha}\cdot h\|^2_{{\rm H}^2}=|c|^2+|d|^2=\|h\|^2_{{\rm H}^2}$ for $h(z)=d+cz^k$.
Since multiplication by $z$ preserves the ${\rm H}^2$-norm, it follows that $\|{\bf b}_{\alpha}\cdot h\|_{{\rm 
H}^2}=\|h\|_{{\rm H}^2}$ holds for all $h$ of the form $h(z)=dz^\ell+cz^k$ and therefore, for all $h\in\bH[z]$.
By the limit argument, the equality holds for every $h\in {\rm H}^2$ and therefore ${\bf b}_{\alpha}$ 
is 
left-inner. Then ${\bf b}_{\overline\alpha}={\bf b}^\sharp_{\alpha}$ is also left-inner, so that 
${\bf b}_{\alpha}$ is right-inner and therefore, bi-inner.\qed

\smallskip

It is clear that the power series (following \cite{acs}, we will call it {\em a finite Blaschke product})
\begin{equation}
B(z)={\bf b}_{\alpha_1}(z){\bf b}_{\alpha_2}(z)\cdots {\bf b}_{\alpha_m}(z)\phi,\qquad 
(\alpha_1,\ldots,\alpha_m\in\mathbb B, \;  |\phi|=1)
\label{9.6}   
\end{equation}
is bi-inner. Its zeros are contained in the union of conjugacy classes $[\alpha_1],\ldots,[\alpha_m]$
and we can use Algorithm \ref{A:4.1} to construct the  spherical divisors of $B$ corresponding to 
each class.
Now we will consider the question related to Algorithm \ref{A:8.1}: {\em to construct a bi-inner power series with 
the prescribed (say, left) zero structure}. By Remark \ref{R:7.10}, in case only finitely many conjugacy
classes are involved, the latter question reduces to the following: {\em given a polynomial 
\begin{equation}
G=\bp_{\alpha_1}\bp_{\alpha_2}\ldots \bp_{\alpha_m},\qquad \alpha_1,\ldots,\alpha_m\in\mathbb B,
\label{9.7}
\end{equation}
find a power series $H\in{\rm H}^2$ with no zeros in $\mathbb B$ so that $f=G\cdot H$ be inner}.

\medskip
\noindent
In case $\alpha_i\alpha_j=\alpha_j\alpha_i$ for $i,j=1,\ldots,m$
(i.e., $\alpha_1,\ldots,\alpha_m$ belong to the same two-dimensional subspace of $\bH$),
we can take $H:={\bf k}_{\overline{\alpha}_1}{\bf k}_{\overline{\alpha}_2}\ldots {\bf 
k}_{\overline{\alpha}_m}$. Some reduction is possible if $G$ has a real zero $x$ 
(i.e., $G=\bp_x\cdot \widetilde G$) or a  spherical zero $[\alpha]$ ($G=\cX_{[\alpha]}\cdot\widetilde G$).
In both cases, it suffices to find $\widetilde H$ so that $\widetilde G\cdot \widetilde H$ is inner
and then let $H={\bf k}_x\cdot \widetilde H$ (in the first case) or $H={\bf k}_{\alpha}\cdot {\bf 
k}_{\overline\alpha}\cdot \widetilde H$ (in the second case). It thus turns out that the general problem reduces 
to the one where $G$ does not have real or spherical zeros. The case where $\alpha_1,\ldots,\alpha_m$ are
pairwise non-equivalent (i.e., has $m$ simple roots) has been handled in \cite{acs}. The next theorem 
settles the general case. 
\begin{theorem}
Given $G$ as in \eqref{9.7}, there exist $\beta_1,\ldots,\beta_m$ ($\beta_i\in[\alpha_i]$) so that 
the power series $G\cdot {\bf k}_{\beta_1}\cdot {\bf k}_{\beta_2}\ldots{\bf k}_{\beta_m}$ is a finite Blaschke 
product.
\label{T:9.7}
\end{theorem}	
{\bf Proof:} We will assume (without loss of generality) that $\alpha_1,\ldots,\alpha_m$ are all 
non-real and prove the statement by induction. For $m=1$, we let $\beta_1=\overline{\alpha}_1$
and the statement follows by Proposition \ref{P:8.3}. Assume that we have found 
$\gamma_i, \delta_i\in[\alpha_i]$ ($i=1,\ldots,m-1$) so that 
\begin{equation}
\bp_{\alpha_1}\cdot \bp_{\alpha_2}\cdots \bp_{\alpha_{m-1}}\cdot {\bf k}_{\delta_1}
\cdot {\bf k}_{\delta_2}\cdots{\bf k}_{\delta_{m-1}}
={\bf b}_{\gamma_1}\cdot {\bf b}_{\gamma_2}\cdots {\bf b}_{\gamma_{m-1}}\varphi.
\label{9.8}
\end{equation}
Since evaluation functional at real points is multiplicative, evaluating both sides in \eqref{9.8} 
at zero gives
$$
\varphi=(\gamma_1\gamma_2\cdots\gamma_{m-1})^{-1}(\alpha_1\alpha_2\cdots\alpha_{m-1}),
$$
and since $\gamma_i\in[\alpha_i]$, it follows that $|\varphi|=1$. 
Given equality \eqref{9.8} and  given $\alpha_m$, we want
to find $\beta_i\in[\alpha_i]$ ($i=1,\ldots,m$) and $\gamma_m\in[\alpha_m]$ such that 
\begin{equation}
\bp_{\alpha_1}\cdot \bp_{\alpha_2}\cdots \bp_{\alpha_{m-1}}\cdot \bp_{\alpha_m}
\cdot {\bf k}_{\beta_1}
\cdot {\bf k}_{\beta_2}\cdots{\bf k}_{\beta_{m}}
={\bf b}_{\gamma_1}\cdot {\bf b}_{\gamma_2}\cdots {\bf b}_{\gamma_{m-1}}\cdot{\bf b}_{\gamma_{m}}
\psi,
\label{9.9}
\end{equation}
where the unimodular factor $\psi$ is given by $\psi=\gamma_{m}^{-1}\varphi \, \alpha_m$. 
Multiplying both parts in \eqref{9.8} by $\varphi^{-1}{\bf b}_{\gamma_{m}}\psi$ on the 
right and  comparing the resulting equality with \eqref{9.9}, we see that \eqref{9.9} is equivalent 
to
$$
\bp_{\alpha_1}\cdots \bp_{\alpha_{m-1}}\cdot {\bf k}_{\delta_1}
\cdots{\bf k}_{\delta_{m-1}}\varphi^{-1}{\bf b}_{\gamma_{m}}\psi=
\bp_{\alpha_1}\cdots \bp_{\alpha_{m-1}}\cdot \bp_{\alpha_m}
\cdot {\bf k}_{\beta_1}\cdots{\bf k}_{\beta_{m}},
$$
which in turn, is equivalent to
\begin{equation}
{\bf k}_{\delta_1}\cdots{\bf k}_{\delta_{m-1}}\varphi^{-1}{\bf b}_{\gamma_{m}}\psi=
\bp_{\alpha_m}\cdot {\bf k}_{\beta_1}\cdots{\bf k}_{\beta_{m}}.
\label{9.10}
\end{equation}
The construction of $\beta_1,\ldots,\beta_m$, $\gamma_m$ and $\psi$ subject to identity \eqref{9.10}
is as follows: letting $\alpha^\prime_1=\alpha_m$ and $\phi_1=1$, we compute 
$\alpha^\prime_k$ and $\phi_k$ by the recursive formulas
\begin{align}
\alpha^\prime_{k+1}&=(1-\overline{\alpha^\prime_k} \, \phi_k^{-1}\delta_{k}\phi_k)\cdot 
\alpha^\prime_k\cdot (1-\overline{\alpha^\prime_k} \, \phi_k^{-1}\delta_{k}\phi_k)^{-1},
\label{9.11}\\
\phi_{k+1}&=\phi_k\cdot (1-\phi_k^{-1}\delta_{k}\phi_k 
\overline{\alpha^\prime_k})(1-\overline{\alpha^\prime_k}\phi_k^{-1}\delta_{k}\phi_k)^{-1},
\label{9.12}
\end{align}
for $k=2,\ldots,m-1$. Then we let 
\begin{align}
\beta_k&=(1-\overline{\alpha^\prime_k} \, \phi_k^{-1}\delta_{k}\phi_k)\cdot
\phi_k^{-1}\delta_{k}\phi_k\cdot (1-\overline{\alpha^\prime_k} \, \phi_k^{-1}\delta_{k}\phi_k)^{-1}
\quad (k=1,\ldots,m-1),\label{9.13}\\
\beta_m&=\overline{\alpha^\prime_m},\quad \psi=\varphi\phi_m,\quad \gamma_m=\psi 
\alpha^\prime_m\psi^{-1}.
\label{9.14}  
\end{align}
It follows from \eqref{9.11} that $\alpha^\prime_{k+1}\sim\alpha^\prime_{k}$ for all $k=1,\ldots,m-1$
and in particular, $\alpha^\prime_{m}\sim\alpha^\prime_{1}=\alpha_m$. By \eqref{9.14}, we also have
$\beta_m,\gamma_m\in [\alpha_m]$. From \eqref{9.13} and the induction hypothesis we conclude that 
$\beta_k\sim\delta_k\in[\alpha_k]$ for $k=1,\ldots,m-1$. By \eqref{9.12}, $|\phi_{k+1}|=|\phi_k|$
for $k=1,\ldots,m-1$ and therefore, $|\phi_{m}|=1$. We also conclude from \eqref{9.13} that 
\begin{equation}
\phi_{m}\bp_{\alpha^\prime_m}\cdot{\bf k}_{\beta_{m}}=
\phi_{m}{\bf b}_{\alpha^\prime_m}=\varphi^{-1}\psi {\bf b}_{\alpha^\prime_m}=
\varphi^{-1}{\bf b}_{\gamma_m}\psi.
\label{9.15}
\end{equation}
To show that \eqref{9.10} holds we first verify equalities
\begin{equation}
\phi_k\alpha^\prime_k=\phi_{k+1}\alpha^\prime_{k+1},\quad
\delta_{k}\phi_k=\phi_{k+1}\beta_k,\quad \phi_k+\delta_{k}\phi_k \alpha^\prime_k=
\phi_{k+1}(1+\alpha^\prime_{k+1}\beta_k)
\label{9.16}    
\end{equation}  
for $k=1,\ldots,m-1$. Indeed, letting for short $\Delta_k=1-\overline{\alpha^\prime_k} \, 
\phi_k^{-1}\delta_{k}\phi_k$, we have from \eqref{9.11}--\eqref{9.13},
\begin{align*}
\phi_{k+1}\alpha^\prime_{k+1}&=\phi_k\cdot (1-\phi_k^{-1}\delta_{k}\phi_k
\overline{\alpha^\prime_k})\alpha^\prime_k\Delta_k^{-1}=
\phi_k \alpha^\prime_k \Delta_k\Delta_k^{-1}=\phi_k \alpha^\prime_k,\\
\phi_{k+1}\beta_k&=\phi_k\cdot (1-\phi_k^{-1}\delta_{k}\phi_k
\overline{\alpha^\prime_k})\phi_k^{-1}\delta_{k}\phi_k\Delta_k^{-1}=
\phi_k \phi_k^{-1}\delta_{k}\phi_k\Delta_k\Delta_k^{-1}=\delta_{k}\phi_k,\\
\phi_{k+1}(1+\alpha^\prime_{k+1}\beta_k)&=
\phi_k\cdot (1-\phi_k^{-1}\delta_{k}\phi_k
\overline{\alpha^\prime_k})(\Delta_k^{-1}(1+\Delta_k 
\alpha^\prime_{k}\phi_k^{-1}\delta_{k}\phi_k \Delta_k^{-1})\\
&=(\phi_k-\delta_{k}\phi_k
\overline{\alpha^\prime_k})(1+\alpha^\prime_{k}\phi_k^{-1}\delta_{k}\phi_k)\Delta_k^{-1}\\
&=\phi_k-\delta_{k}\phi_k \overline{\alpha^\prime_k}+(\phi_k-\delta_{k}\phi_k
\overline{\alpha^\prime_k})(\alpha^\prime_k+\overline{\alpha^\prime_k})
\phi_k^{-1}\delta_{k}\phi_k\Delta_k^{-1}\\
&=\phi_k-\delta_{k}\phi_k 
\overline{\alpha^\prime_k}+(\alpha^\prime_k+\overline{\alpha^\prime_k})
\delta_{k}\phi_k\Delta_k\Delta_k^{-1}\\
&=\phi_k-\delta_{k}\phi_k\overline{\alpha^\prime_k}+\delta_{k}\phi_k
(\alpha^\prime_k+\overline{\alpha^\prime_k})=\phi_k+\delta_{k}\phi_k \alpha^\prime_k.
\end{align*}
The consequence of relations \eqref{9.16} is the polynomial identity
$$
(1-z\delta_{k})\phi_k(z-\alpha^\prime_k)=\phi_{k+1}(z-\alpha^\prime_{k+1})(1-z\beta_k),
$$
which being multiplied by ${\bf k}_{\delta_k}(z)$ on the left and by ${\bf k}_{\beta_k}(z)$
on the right, implied 
$$
\phi_k\bp_{\alpha^\prime_k}{\bf k}_{\beta_k}={\bf k}_{\delta_k}\phi_{k+1}\quad\mbox{for}\quad 
k=1,\ldots,m-1.
$$
Combining the latter equalities with \eqref{9.15} gives
\begin{align*}
\bp_{\alpha_m}\cdot {\bf k}_{\beta_1}\cdot{\bf k}_{\beta_2}\cdot{\bf k}_{\beta_3}\cdots{\bf 
k}_{\beta_{m}}&=
\phi_1\bp_{\alpha^\prime_1}\cdot{\bf k}_{\beta_1}\cdot{\bf k}_{\beta_2}\cdot{\bf 
k}_{\beta_3}\cdots{\bf k}_{\beta_{m}}\\
&={\bf k}_{\delta_1}\phi_2 \bp_{\alpha^\prime_2}\cdot{\bf k}_{\beta_2}\cdot{\bf k}_{\beta_3}\cdots{\bf 
k}_{\beta_{m}}\\
&={\bf k}_{\delta_1}\cdot{\bf k}_{\delta_2}\phi_3 \bp_{\alpha^\prime_3}\cdot{\bf 
k}_{\beta_3}\cdots{\bf k}_{\beta_{m}}=\ldots \\
&={\bf k}_{\delta_1}\cdot{\bf k}_{\delta_2}\cdots{\bf 
k}_{\delta_{m-1}}\phi_{m}\bp_{\alpha^\prime_m}\cdot{\bf k}_{\beta_{m}}\\
&={\bf k}_{\delta_1}\cdots{\bf k}_{\delta_{m-1}}\varphi^{-1}{\bf b}_{\gamma_{m}}\psi
\end{align*}
which confirms \eqref{9.10} and completes the proof of the theorem.\qed

\smallskip

\bibliographystyle{amsplain}

\end{document}